\documentclass{report}

\usepackage{amsfonts,latexsym}
\usepackage{epsfig}
\usepackage{amsmath}
\usepackage{amssymb}
\usepackage{graphicx}
\usepackage{caption}
\usepackage{subfig}
\usepackage{sidecap}
\usepackage{caption}

\usepackage[latin1]{inputenc}
\usepackage[english]{babel}
\newtheorem{theorem}{Theorem}[section]
\newtheorem{lemma}{Lemma}[section]

\newtheorem{definition}{Definition}

\newtheorem{assumption}{Assumption}[section]

\begin{document}

\title{Error Analysis  \\  of Domain Decomposition method for\\ 4D Variational Data Assimilation (4DVAR DA): Consistence, Convergence and Stability.}
\author{Rosalba Cacciapuoti, Luisa D'Amore}
\maketitle
\date{}
\noindent Abstract: We prove consistence,  convergence and stability  of the Domain Decomposition in space and time method for 4DVAR Data Assimilation problem (DD - 4D VAR). We give the condition number of the DD - 4D VAR method and  experimentally validate the analysis on the Shallow Water Equations.
\section{4DVAR DA formulation  }

If $\Omega \subset \mathbb{R}^{n}$, $n\in \mathbb{N}$,  is a spatial  domain with a Lipschitz boundary, let:
\begin{equation}\label{modelloDA}
\left\{ \begin{array}{ll}
u^{\mathcal{M}}(t+h,x)=\mathcal{M}_{t,t+h}[u(t,x)] & \textrm{$\forall x \in \Omega$, $t,t+h \in [0,T]$, $(h >0)$} \\
u^{\mathcal{M}}(t_{0},x)=u_{0}(x) & \textrm{$ t_{0}\equiv 0, \ \ x\in \Omega$}\\
u^{\mathcal{M}}(t,x)=f(x) & \textrm{$ x\in \partial \Omega$, $\forall t \in [0,T]$}\\
\end{array}, \right.
\end{equation}
be a symbolic description of the  model of interest where
\begin{equation}\label{fun_u}
u^{\mathcal{M}}:(t,x) \in [0,T] \times \Omega \mapsto u^{\mathcal{M}}(t,x)= [u^{\mathcal{M}}[1](t, x), u^{\mathcal{M}}[2](t, x),\ldots, u^{\mathcal{M}}[pv](t, x)],
\end{equation}
\noindent 
is the state function of $\mathcal{M}$ with $pv\in \mathbb{N}$ the number of physical variables, $f$ is a known function defined on the boundary $ \partial \Omega$;let
$$v:(t,x) \in [0,T] \times \Omega \mapsto v(t,x),$$
be the observations function and
$$\mathcal{H}: u^{\mathcal{M}}(t,x) \mapsto v(t,x), \ \ \ \ \ \forall (t,x) \in [0,T] \times \Omega,$$
denote the non-linear observations mapping. To simplify future treatments we assume $pv \equiv 1$. \\

\begin{definition}\label{discreto_dominio}(Discretization of  $\Omega\times \Delta$)
Let 
\begin{equation*}
 \Omega_I\equiv   \{x_{\tilde{i}}\}_{\tilde{i}\in I}\subset \Omega
\end{equation*}
be the discretization of  $\Omega$ where 
\begin{equation}\label{setI}
I=\{1,\ldots,N_p\} \ \textrm{and} \ N_{p}=|I| \footnote{We refer to $|I|$ as cardinality of set $I$.}
\end{equation}
are respectively the set of indices of nodes in $\Omega$ and its cardinality i.e. the number of inner nodes in $\Omega$. Let
\begin{equation*}
  \Delta_{K}  \equiv    \{t_{\tilde{k}}\}_{\tilde{k}\in K} \subset \Delta
\end{equation*}
be the  discretization of  $\Delta$ where
\begin{equation}\label{K}
    K=\{1,\ldots,N\} \ \textrm{and} \ N=|K| 
\end{equation}
are respectively the set of indices of the time variable in $\Delta$ and its cardinality i.e. number of time instants in  $\Delta$. Consequently, we refer to
\begin{equation}\label{tot_discretization}
 \Omega_{I}\times \Delta_K\equiv    \{(x_{\tilde{i}},t_{\tilde{k}})\}_{\tilde{i}\in I; \ \tilde{k} \in K}\subset \Omega\times \Delta
\end{equation}
as the discrete domain.\\
\end{definition}

\noindent We introduce 4D-DA variational formulation, i.e. the 4D-VAR DA problem.
\begin{definition}(The 4DVAR DA problem). Let $\Omega\subset \mathbb{R}^n$ and $\Delta\subset \mathbb{R}$ be spatial domain and time interval. The 4DVAR DA  problem concerns the computation of:
\begin{equation}\label{varDA}
\mathbf{u}^{DA}=argmin_{u\in \mathbb{R}^{N_{p}\cdot N}}J(u),
\end{equation}
solution of 4DVar problem on $\Omega\times \Delta$ with
\begin{equation}\label{funzionale}
J(u)=\alpha \|u-u^{M}\|_{\textbf{B}^{-1}}^{2}+\|Gu-y\|_{\textbf{R}^{-1}}^{2},
\end{equation}
where
\begin{itemize}

\item $N_{p}$: is the number of nodes in $\Omega \subset \mathbb{R}^{n}$ defined in (\ref{setI});
\item $n_{obs}$: is the number of observations in $\Omega$; where $n_{obs} <<N_{p}$, 
\item $N$: is the number of time instants in $\Delta$ defined in (\ref{K});
\item $\alpha$: is the regularization parameter
\item 
$u_{0}=\{u_{0,j}\}_{j=1,\ldots,N_{p}}\equiv \{u_{0}(x_{j})\}_{j=1,\ldots,N_{p}} \in \mathbb{R}^{N_{p}}$:
 is the state at time $t_{0}$;
\item the operator 
\begin{equation}\label{model_disc}
M_{l-1,l}\in \mathbb{R}^{N_{p} \times N_{p}}, \ \ \ l=1,\ldots,N,
\end{equation}
is the discretization of the linear approximation of $\mathcal{M}_{t_{l-1},t_{l}}$ from $t_{l-1}$ to $t_{l}$; 
\item the operator 
\begin{equation}\label{model_disc_tot}
M\in \mathbb{R}^{N_{p} \times N_{p}}
\end{equation}
is the discretization of the linear approximation of $\mathcal{M}$ from $t_{0}$ to $t_{N}$;
\item the matrix
\begin{equation}\label{background}
u^{M}:=\{u_{j,l}^{{M}}\}_{j=1,\ldots,N_{p};l=1,\ldots,N} \equiv \{u^{M}(x_{j},t_{l})\}_{j=1,\ldots,N_{p};l=0,1,\ldots,N-1} \in \mathbb{R}^{N_{p}\times N}
\end{equation}
is the solution of  $M$  i.e. the background;
\item 
$y:= \{y(z_{j},t_l)\}_{j=1,\ldots,n_{obs};l=0,1,\ldots,N-1}\in \mathbb{R}^{ n_{obs}\times N}$:
are the observations; 

\item 
$H_{l}\in \mathbb{R}^{n_{obs} \times N_{p}}, \ \ \ l=0,\ldots,N-1$
: is the linear approximation of observation mapping $\mathcal{H}$;
\item $G\equiv G_{N-1}\in \mathbb{R}^{(N \times n_{obs})\times N_{p} }$:
\begin{displaymath}\label{mat_G}
G_{l}=\left\{ \begin{array}{ll}
\left[\begin{array}{ll}
H_{0}\\ H_{1}\\ \vdots \\H_{l-1} \end{array}\right] & \textrm{$l>1$} \\
\\
H_{0}& \textrm{$ l=1$}
\end{array}, \right.
\end{displaymath}
\item $\textbf{R}=diag(\textbf{R}_{0},\textbf{R}_{1},\ldots,\textbf{R}_{N-1})$  and \textbf{B}$=\textbf{V}\textbf{V}^{T}$: are covariance matrices of the errors on  observations and  background, respectively.
\end{itemize}

\end{definition}
\section{The space and time Domain Decomposition method: DD--4DVAR }\label{section_DD_4DVAR}

\begin{figure}
\centering
{\includegraphics[width=0.7\textwidth]{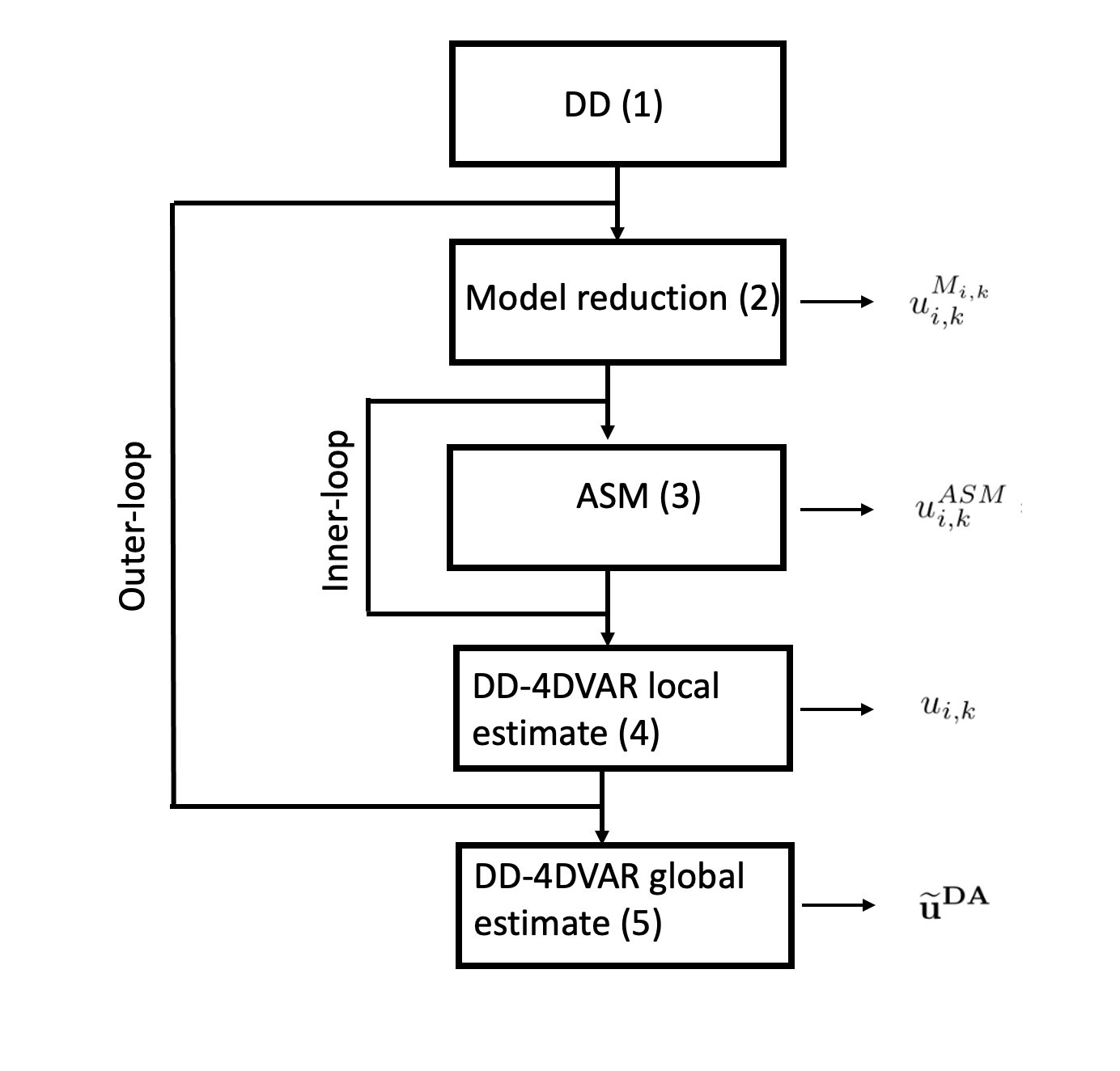}}\\
\caption{Schematic description of DD--4DVAR algorithm. The modules of DD--4DVAR: DD, Model Reduction, ASM, DD--4DVAR local solution and DD--4DVAR global solution are identified and the Arabic numbers in parentheses refer to the corresponding module in section \ref{section_DD_4DVAR}. For each module we report the corresponding solution.}
\label{schema1}
\end{figure}

The  DD method consists in decomposing the domain of computation $\Omega \times\Delta$ into subdomains   and solving reduced forecast models and local 4DVAR DA problems in \cite{dd-4dvar}. 
The modules of DD method  are: domain decomposition of $\Omega\times \Delta$, Model Reduction, ASM based on Additive Schwarz method \cite{ASM}, DD-4DVAR local solution and DD-4DVAR global solution. The schematic description of DD-4DVAR algorithm is reported in Figure \ref{schema1}.
\\

\begin{enumerate}

\item \textbf{Domain Decomposition  of $\Omega\times \Delta$.}\\ 
We describe DD of $\Omega\times \Delta$ and $\Omega_I\times \Delta_k$.
\begin{itemize}
    \item
DD of $\Omega\times \Delta$ consists of:
\\
decomposition of  $\Omega\subset \mathbb{R}^n$ into a sequence of subdomains $\Omega_{i}$ such that:
\begin{equation}\label{DD_space}
\Omega=\bigcup_{i=1}^{N_{sub}}\Omega_{i}
\end{equation}
and definition of set of indices of  subdomains  adjacent to $\Omega_i$ and its cardinality, as follows 
\begin{equation}\label{ni}
J_i\subset \{1,\ldots,N_{sub}\} \quad \textrm{and} \quad ad_i=|J_i|
\end{equation}
in particular, $ad_i$ is the number of  subdomains adjacent to $\Omega_i$.
\\\
For $i=1,\ldots,N_{sub}$, definition of overlap regions $\Omega_{ij}$ as follows 
\begin{equation}\label{overlap_region}
\Omega_{ij}:=\Omega_{i}\cap \Omega_{j}\neq \emptyset \quad \textrm{$\forall j\in J_i$}.
\end{equation}
 Definition of interfaces  $\Omega_{i}$, for $i=1,\ldots,N_{sub}$:
\begin{equation}\label{interfacce}
\Gamma_{ij}:=\partial \Omega_{i} \cap \Omega_{j} \quad \textrm{ $ j\in J_i$}.
\end{equation}
 Decomposition of time interval $\Delta\subset \mathbb{R}$ into a sequence of time interval $\Delta_k$ such that:
$$\Delta=\bigcup_{k=1}^{N_{t}}\Delta_{k}.$$
Consequently. we define
\begin{equation}\label{local_domain}
\{\Omega_i\times \Delta_k \}_{i=1,\ldots,N_{sub};k=1,\ldots,N_t}
\end{equation} as local domains.
\item DD of $\Omega_I\times \Delta_K$ defined in (\ref{tot_discretization}) consists of:
\\
identification of inner nodes of subdomains  $\{\Omega_i\}_{i=1,\ldots,N_{sub}}$: for $i=1,\ldots,N_{sub}$ 
\begin{equation*}
        \{x_{\tilde{i}}\}_{\tilde{i}\in I_i}\subset \Omega_i
\end{equation*}
are inner nodes of $\Omega_i$ where $I_i$ is set defined as follows
$$I_i:=\left\{(i-1)\times \frac{N_p}{N_{sub}}+1,\ldots, i\times \frac{N_p}{N_{sub}}+\frac{\delta}{2}\right\}$$ 
such that 
\begin{equation}\label{insieme_I}
    I=\bigcup_{i=1}^{N_{sub}} I_i.
\end{equation}
with $I$ set of indices of inner nodes in $\Omega$ defined in (\ref{setI}).
\\

\noindent Identification of inner nodes of overlap regions $\{\Omega_{ij}\}_{i=1,\ldots,N_{sub},j\in J_i}$: for $i=1,\ldots,N_{sub}$   
\begin{equation*}
  \Omega_{I_i}\equiv       \{x_{\tilde{j}}\}_{\tilde{j}\in I_{ij}}\subset \Omega_{ij} \quad \forall j\in J_i
\end{equation*}
are inner nodes of $\Omega_{ij}$ where 
\begin{equation}\label{Iij}
  I_{ij}:=I_i\cap I_j \neq \emptyset \quad \forall j \in J_i
\end{equation}
and 
\begin{equation}\label{delta}
    \delta:=|I_{i,j}|
\end{equation}
are respectively the set of indices of nodes in overlap region and its  cardinality i.e. the number of inner nodes in overlap region $\Omega_{ij}$ defined in (\ref{overlap_region}).
Consequently, for $i=1,\ldots,N_{sub}$ we define the cardinality of $I_i$ as follows:
\begin{equation}\label{nloc1}
N_{loc}:=|I_i|=\frac{N_p}{N_{sub}}+\frac{\delta}{2}
\end{equation}
i.e. number of inner nodes of spatial subdomain $\Omega_i$.
\\
\\
Identification of time instant in  $\{\Delta_k\}_{k=1,\ldots,N_{t}}$: for $k=1,\ldots,N_{t}$ 
\begin{equation*}
  \Delta_{K_k}\equiv      \{t_{\tilde{k}}\}_{\tilde{k}\in K_k}\subset \Delta.
\end{equation*}
are the time instant  in $\Delta_k$ where $K_k$ is  defined as follows:  
\begin{equation*}
    K_{k}:=\left\{(k-1)\times \frac{N}{N_{t}},\ldots, k\times \frac{N}{N_{t}}\right\}
    \end{equation*}
   where
    \begin{equation*}
         K_k\cap K_{k+1}=\left\{k\times \frac{N}{N_{t}}\right\}\neq  \emptyset \quad \forall k=1,\ldots,N_t-1
    \end{equation*}
    and 
\begin{equation*}
    N_k:=|K_k|=\frac{N}{N_t}
\end{equation*}
is  cardinality of $K_k$ i.e. the number of time instants in $K_k$, such that 
\begin{equation}\label{insieme_K}
    K=\bigcup_{k=1}^{N_{t}} K_k.
\end{equation}
with $K$ set defined in (\ref{K}).



Consistently with Definition \ref{discreto_dominio}, for $i=1,\ldots,N_{sub}$, $k=1,\dots,N_t$
\begin{equation*}
   \Omega_{I_i}\times\Delta_{K_k}:= \{(x_{\tilde{i}},t_{\tilde{k}})\}_{\tilde{i}\in I_{i}; \tilde{k}\in K_k}\subset \Omega_i\times \Delta_k
\end{equation*}
is local discrete domain.

\end{itemize}
In order to compute local estimates on local domains defined in (\ref{local_domain}) we introduce the following notation.
For $i=1,\ldots,N_{sub}$, $k=1,\dots,N_t$, we pose
\begin{equation*}
    z_{i,k}=\{z({\tilde{i}},{\tilde{k}})\}_{\tilde{i}\in I_i; \ \tilde{k} \in K_k}\in \mathbb{R}^{N_{loc}\times N_k}
\end{equation*}
i.e. a vector defined on local domain $\Omega_i\times \Delta_k$.
\\

We are able to define the restriction and extension operator.
\begin{definition}\label{def_rest}(Restriction Operator) 
For $i=1,\ldots,N_{sub}$ we define restriction matrices $R_{i}$, $R_{ij}$ on $\Omega_{i}$, $\Omega_{ij}$, respectively. 
Given $x\in \mathbb{R}^{N_p\times N}$ and $z\in \mathbb{R}^{N_p\times 1}$ we define restriction of $x$ to $\Omega_i$ :
\begin{equation}\label{matrice_R}
 x/ \Omega_{{i}}:= R_{i} x=\{x({\tilde{i}},\tilde{k})\}_{\tilde{i}\in {I}_{i},\tilde{k}\in K_k} \in \mathbb{R}^{N_{loc}\times 1} 
 \end{equation}
 $$ 
x/ \Omega_{{ij}}:= R_{ij} x=
\{x({\tilde{j}},\tilde{k})\}_{\tilde{j}\in {I}_{ij},\tilde{k}\in K_k} \in \mathbb{R}^{\delta\times 1}
$$
and restriction of $y$ to $\Omega_i\times \Delta_k$
$$ z/ (\Omega_{{i}}\times \Delta_k):=\{z(\cdot,\tilde{k})\}_{\tilde{k}\in K_k} / \Omega_{{i}}= R_{i}\cdot \{z(\cdot,\tilde{k})\}_{\tilde{k}\in K_k}=\{z({\tilde{i}},\tilde{k})\}_{\tilde{i}\in {I}_{i},\tilde{k}\in K_k} \in \mathbb{R}^{N_{loc}\times N_k}$$ $$ z/ (\Omega_{{ij}}\times \Delta_k):=\{z(\cdot,\tilde{k})\}_{\tilde{k}\in K_k} / \Omega_{{ij}}= R_{ij}\cdot \{z(\cdot,\tilde{k})\}_{\tilde{k}\in K_k}=\{z({\tilde{i}},\tilde{k})\}_{\tilde{i}\in {I}_{ij},\tilde{k}\in K_k} \in \mathbb{R}^{N_{loc}\times N_k}$$
where $I_i$ and $I_{ij}$ are respectively set of indices of inner nodes in $\Omega_i$ and $\Omega_{ij}$, $\forall j\in J_i$.
\end{definition}

\begin{definition}\label{def_ex}(Extension operator) 
We define the Extension Operator  (EO). If  $x\in \mathbb{R}^{N_{loc}\times N_k}$, it is
$$ EO(x):=R_{i}^Tx= \left \{
\begin{array}{cc}
  x({\tilde{i}},\tilde{k}) & \textit{if}\ ({\tilde{i}},\tilde{k})\in I_i\times K_k  \\
  0 & elsewhere
\end{array}
\right .
$$
where $R_{i}^T$ is the transpose of $R_i$ in (\ref{matrice_R}) and
$ EO(x) \equiv {x^{EO}}$.\\
\end{definition}
At this point, we introduce the outer-loop solving the nonlinear aspects of the assimilation
problem which for DD-4DVAR, in general for 4DVAR techniques, includes the integration of the nonlinear model. \\
For $n=0,1,\ldots, n_{stop}$, do\\
\item \textbf{Model Reduction:} for $ i=1,\ldots,N_{sub}$; $k=1,\ldots,N_{t}$, posed $u_{i,k}^{0}:= \{u^{{M}}(x_{\tilde{i}},t_{\tilde{k}})\}_{\tilde{i}\in {I}_{ij},\tilde{k}\in K_k}$ i.e. by using the background as local initial values, let $u_{i,k}^{{M}_{i,k},n+1}$ be the solution of the local model $(P_{i,k}^{{M}_{i,k},n})_{i=1\ldots,N_{sub},k=1,\ldots,N_{t}}$: 
\begin{equation}\label{co}
(P_{i,k}^{{M}_{i,k},n})_{i=1\ldots,N_{sub},k=1,\ldots,N_{t}}\ : \quad 
\left\{\begin{array}{ll}
u_{i,k}^{{M}_{i,k},n}=M_{i,k}\cdot u_{i,k-1}^{n}+b_{i,k}^{n},\\
u_{i,k-1}^{n}=u_{i,k-1}^{{M}_{i,k},n}\\
u_{i,k}^{n}/\Gamma_{ij}=u_{j,k}^{n}/\Gamma_{ij}, \quad j\in J_i
\end{array}\right.
\end{equation}
where $u_{i,k}^{M}$, $b_{i,k}^{n}$ are respectively the background in $\Omega_i \times \Delta_{k}$, the vector accounting 
boundary conditions of $\Omega_{i}$ and
\begin{equation}\label{local_matrix}
    M_{i,k}:=M_{k}/\Omega_i
\end{equation}
the restriction to $\Omega_i$ of the matrix 
\begin{equation}\label{matriceM}
M_{k}\equiv M_{\bar{s}_{k-1},\bar{s}_{k}}:= M_{\bar{s}_{k-1},\bar{s}_{k-1}+1}\cdots M_{\bar{s}_{k}-1,\bar{s}_{k}}.
\end{equation}
where 
\begin{equation}\label{initial_indices}
    \bar{s}_{k}:=\sum_{j=1}^{k-1}N_{j}-(N-1) \quad \textrm{and}\quad \bar{s}_{0}:=0
    \end{equation}
are respectively the first index of $\Delta_{K_K}$ and $\Delta_{K_1}$ and $b_{i,k}^{n}$ is vector computed from information on boundary of $\Omega_i$ at  $\Delta_k$ that depends on the discretization scheme used.
\\

Let:
\begin{equation}\label{co2}
(P_{i,k}^{n})_{i=1,\ldots,N_{sub}\,,k=1,\ldots,N_{t}} \ : 
u_{i,k}^{ASM,n}=\arg\min_{u_{i,k}^{n}}\,\,\textbf{J}_{i,k}(u_{i,k}^{n}) 
\end{equation}
be the local  4DVAR DA model with 
\begin{equation}\label{local_operator}
\mathbf{J}_{i,k}(u_{i,k}^{n})=\textbf{J}(u_{i,k}^{n})/(\Omega_{i}\times \Delta_k)+\mathcal{O}_{ij}.
\end{equation}
We let
\begin{equation}\label{J_ris0}
\mathcal{O}_{ij}=\sum_{ j\in J_i}\beta_{j}\cdot \|u_{i,k}^{n}/\Omega_{ij}-u_{j,k}^{n}/\Omega_{ij}\|_{\textbf{B}_{ij}^{-1}}^{2}
\end{equation}
be the overlapping operator on $\Omega_{ij}$, and  

\begin{equation*}
\textbf{J}_{i,k}(u_{i,k}^n)/(\Omega_{i}\times \Delta_k)=  \alpha_{i,k} \cdot \|u_{i,k}^{n}-u_{i,k}^{{M}_{i,k},n}\|_{\textbf{B}_{i}^{-1}}+\|{G}_{i,k}u_{i,k}^{n}-y_{i,k}\|_{\textbf{R}_{i}^{-1}}^{2}
\end{equation*}
be the overlapping operator on $\Gamma_{ij}$, and  

\begin{equation*}
\textbf{J}_{i,k}(u_{i,k}^n)/(\Omega_{i}\times \Delta_k)=  \alpha_{i,k} \cdot \|u_{i,k}^{n}-u_{i,k}^{{M}_{i,k},n}\|_{\textbf{B}_{i}^{-1}}+\|{G}_{i,k}u_{i,k}^{n}-y_{i,k}\|_{\textbf{R}_{i}^{-1}}^{2}
\end{equation*}
be the restriction of $\mathbf{J}$ to $\Omega_{i}\times \Delta_k$. Parameters $\alpha_{i,k}$ and  $\beta_{j}$ in (\ref{J_ris0}) denotes the regularization parameters. Following, we let $\alpha_{i,k}=\beta_{j}=1$, $j\in J_i$.\\

\item \textbf{ASM}: according to ASM \cite{ASM} we compute solution of $( P_{i,k}^{n})_{i=1,\ldots,N_{sub}\,,k=1,\ldots,N_{t}}$.  Gradient of $\textbf{J}_{i,k}$ is \cite{JCP,articolo2}:
\begin{equation}\label{gradiente}
\begin{array}{ll}
\nabla \textbf{J}_{i,k}(\textbf{w}_{i,k}^{n})=&(\textbf{V}_{i}^{T}({G}_{i,k})^{T}(\textbf{R}_{i,k})^{-1}{G}_{i,k}\textbf{V}_{i}+ I_{i}+ad_{i}\cdot\textbf{B}_{ij} )\textbf{w}_{i,k}^{n}\\&-\textbf{c}_{i}+\sum_{ j\in J_i}\textbf{B}_{ij} \textbf{w}_{j,k}^{n},
\end{array}
\end{equation}
where 
\begin{equation}\label{w}
\textbf{w}_{i,k}^{n}=\textbf{V}_{i}^{-1}(u_{i,k}^{n}-u_{i,k}^{{M}_{i,k},n}),
\end{equation}
$$\textbf{d}_{i}=(\textbf{v}_{i}-{G}_{i,k}{u}_{i,k}^{{M}_{i,k},n}) \quad  \textbf{c}_{i}=(\textbf{V}_{i}^{T}({G}_{i,k})^{T}(\textbf{R}_{i,k})^{-1}\textbf{d}_{i})$$ 
and $I_{i}$ the identity matrix, $ad_i$ number of  subdomains adjacent to $\Omega_i$ defined in (\ref{ni}). Solution of $( P_{i,k}^{n})_{i=1,\ldots,N_{sub}\,,k=1,\ldots,N_{t}}$ is obtained by requiring $\nabla \textbf{J}_{i,k}(\textbf{w}_{i,k}^{n})=0$. This gives rise to the linear system: 
\begin{equation}\label{system_n}
A_{i,k}\textbf{w}_{i,k}^{n}=c_{i}-\sum_{j\in J_i} \textbf{B}_{ij}\textbf{w}_{j,k}^{n}, 
\end{equation}
where
\begin{equation}\label{Aik}
A_{i,k}=(\textbf{V}_{i}^{T}({G}_{i,k})^{T}\textbf{R}_{i,k}^{-1}{G}_{i,k}\textbf{V}_{i}+ I_{i}+ad_{i}\cdot \textbf{B}_{ij}).
\end{equation}
\\
\noindent (7.1) For each iteration $n$, r.h.s. of (\ref{system_n}) depends on unknown value $\textbf{w}_{j,k}^{n}$ defined on  sub domains $\Omega_{ij}$ adjacent to $\Omega_i$,for $j \in J_i$. Hence, for solving the system in (\ref{system_n}).\\
We introduce an inner-loops related to minimization of local 4DVAR model defined in (\ref{co2}). For $r=0,1,\ldots,\bar{r}$ 
\begin{equation}\label{MPS}
A_{i,k}\textbf{w}_{i,k}^{r+1,n}=c_{i}-\sum_{j \in J_i} \textbf{B}_{ij}\textbf{w}_{j,k}^{r,n}, 
\end{equation}
where at each step $r+1$,  $\Omega_{i}$ \textit{receives}  $\textbf{w}_{i,k}^{r,n}$   from $\Omega_{ij}$, $j \in J_i$ for computing the r.h.s. of (\ref{MPS}) then it  computes  $\textbf{w}_{i,k}^{r+1,n}$ by using  Conjugate Gradient (CG) method and finally  it \textit{sends} $\textbf{w}_{i,k}^{r+1,n}$ to  $\Omega_{ij}$, $j \in J_i$ for updating the  r.h.s. of  (\ref{MPS}) needed to  the next iteration.
$\textbf{w}_{i_j,k}^{0,n}$ is an arbitrary initial value.
\noindent  Finally,  we pose 
$$\textbf{w}_{i,k}^{n}\equiv \textbf{w}_{i,k}^{\bar{r},n},$$
consequently 
\begin{equation}\label{stima_ASM_rstop}
u_{i,k}^{ASM,n}:=u_{i,k}^{ASM,\bar{r}}.
\end{equation}
\item \textbf{DD--4DVAR approximation in $\Omega_i\times\Delta_k$: }final solution update, using (\ref{co}) and (\ref{w}):
\begin{equation}\label{solpara}
\begin{array}{ll}
u_{i,k}^{n+1}
&=u_{i,k}^{{M}_{i,k},n+1}+\textbf{V}_{i}\textbf{w}_{i,k}^{n}= u_{i,k}^{{M}_{i,k},n+1}+[u_{i,k}^{ASM,n}-u_{i,k}^{{M}_{i,k},n}].
\end{array} 
\end{equation}
Endfor $n$ \\

\item \textbf{DD--4DVAR solution in $\Omega\times\Delta$: }computation of DD-4DVAR approximation on $\Omega\times \Delta$: let

\begin{equation}\label{tot_n}
\mathbf{\widetilde{u}^{DD-DA,n}}:= \sum_{i=1}^{N_{sub}}\sum_{k=1}^{N_{t}}({{u}}_{i,k}^{{n}})^{EO}.
\end{equation}
be DD--4DVAR approximation at iteration $n$ where  $({{u}}_{i,k}^{{n}})^{EO}$ is extension to $\Omega\times \Delta$ of DD-4DVAR approximations in $\Omega_i\times \Delta_k$, we define 
\begin{equation}\label{def_uDA}
\mathbf{\widetilde{u}^{DD-DA}}:=\mathbf{\widetilde{u}^{DD-DA,{\bar{n}}}}
\end{equation}
as DD-4DVAR solution in $\Omega\times \Delta$.

\end{enumerate}

\noindent Note that \textbf{B}$_{i}=R_{i}\textbf{B}R_{i}^{T}$ and $ \textbf{B}_{ij}:=\textbf{B}/\Gamma_{ij}=R_{i}\textbf{B}R_{ij}^{T}$ are the restrictions of  covariance matrix $\textbf{B}$, respectively, to  subdomain $\Omega_{i}$ and interface $\Gamma_{ij}$ in (\ref{interfacce}),  ${G}_{i,k}$, $\textbf{R}_{i,k}$ are the restriction of matrices ${G}_{k}:={G}_{\bar{s}_{k}}$
 and ${\textbf{R}_k}:=diag(\textbf{R}_{0},\textbf{R}_{1},\ldots,\textbf{R}_{\bar{s}_{k}})$ to $\Omega_{i}$, and finally  ${u}_{i,k}^{{M}}=R_{i}{u}_{k}^{{M}}$, ${u}_{i,k}^{n+1}/\Gamma_{ij}=R_{ij}{u}_{i,k}^{n+1}$, ${u}_{j,k}^{n}/\Gamma_{ij}=R_{ij}{u}_{k}^{n}$ are the restriction of vectors ${u}_{k}^{b}$, ${u}_{i,k}^{n+1}$, ${u}_{j,k}^{n}$ to $\Omega_{i}$ and $\Gamma_{ij}$,  for $i=1,2,\ldots,N_{sub}$, $j \in J_i$. \\


\noindent  It is worth noting that each local functional  $\textbf{J}_{i,k}$  is obtained starting  from a restriction of the global  functional $\textbf{J}$ and adding a local term defined on the overlapping regions. In addition, regarding the decomposition in time direction, we use   DA  as predictor operator for the local PDE-based model, providing the approximations needed for  solving  the initial value problem on each sub interval. 

\noindent We prove that the minimum of $\mathbf{J}$ can be obtained by patching together all the local solution obtained as minimum of local function $\mathbf{J}_{i,k}$.
\noindent We are able to prove the following result.

\begin{theorem} \label{theorem215}
If $\mathbf{J}$ is convex:

 \begin{equation}\label{thesis}
\mathbf{J}(\mathbf{u}^{DA}) =  \mathbf{J}( \mathbf{\widetilde{u}^{DD-DA}}).
 \end{equation}

\end{theorem}

\noindent \textbf{Proof:} 
\noindent As $ \mathbf{u}^{DA}$ is the global minimum of $\mathbf{J}$ it follows that:

\begin{equation}\label{diseq-minim}
 \mathbf{J}(\mathbf{{u}^{DA}}) \leq  \mathbf{J} \left(\sum_{i=1}^{N_{sub}}\sum_{k=1}^{N_{t}}{(\mathbf{u}}_{i,k}^{\bar{n}})^{EO}\right), \quad \forall \,i,k
\end{equation}

\noindent then, from the (\ref{def_uDA}) it follows that

\begin{equation}\label{diseq-minim2}
 \mathbf{J}(\mathbf{u^{DA}})\leq  \mathbf{J}\left( \mathbf{\widetilde{u}^{DD-DA}}\right )\quad .
\end{equation}

\noindent  Now  we prove that  if $\mathbf{J}$ is convex, then $$\mathbf{J}( \mathbf{u^{DA}} )= \mathbf{J}(\mathbf{\widetilde{u}^{DD-DA}})$$ by reduction to the absurd. Assume that

\begin{equation}\label{PerAbsurd}
\mathbf{J}( \mathbf{u^{DA}} )< \mathbf{J}( \mathbf{\widetilde{u}^{DD-DA}} ). 
\end{equation}

\noindent In particular, 
\begin{equation}\label{zero0}
\mathbf{J}(\mathbf{u^{DA}}) <\mathbf{J}\left( \sum_{i=1}^{N_{sub}}\sum_{k=1}^{N_{t}} \mathbf{\widetilde{u}_{i,k}}\right) \le \sum_{i=1}^{N_{sub}} \sum_{k=1}^{N_{t}} \mathbf{\mathbf{J}_{i,k}}({{u}_{i,k}^{\bar{n}}})  \quad .
\end{equation}
\noindent This means  that

\begin{equation}\label{eq_ROabsurd}
  \mathbf{J}(\mathbf{u^{DA}} ) < \sum_{i=1}^{N_{sub}} \sum_{k=1}^{N_{t}} \mathbf{\mathbf{J}_{i,k}}({{u}_{i,k}^{\bar{n}}}) \quad .
\end{equation}

\noindent From the (\ref{eq_ROabsurd}) and the (\ref{zero0}), it is:

\begin{equation}\label{tre}
   \nabla \mathbf{J}(\mathbf{u}^{DA})<\nabla \sum_{i=1}^{N_{sub}} \sum_{k=1}^{N_{t}} \mathbf{\mathbf{J}_{i,k}}({{u}_{i,k}^{\bar{n}}})=\sum_{i=1}^{N_{sub}} \sum_{k=1}^{N_{t}}\nabla  \mathbf{\mathbf{J}_{i,k}}({{u}_{i,k}^{\bar{n}}}) 
\end{equation}
$\mathbf{u}^{DA}$ is global minimum of $J$ i.e. $\nabla \mathbf{J}(\mathbf{u}^{DA})=0$, from (\ref{tre}) it is 

\begin{equation}\label{Eq_theAbsurd}
\nabla  \mathbf{\mathbf{J}_{i,k}}({u}_{i,k})>0\quad \forall (i,k).
\end{equation}

\noindent The (\ref{Eq_theAbsurd}) is an absurd as the values of ${{u}_{i,k}^{\bar{n}}}$ are the minimum for $\mathbf{\mathbf{J}_{i,k}}$ for $\forall (i,k)$.  Hence, the (\ref{thesis}) is proved.\\
\\
\\

\section{Consistency,  convergence and stability of DD-4DVAR method}\label{conv_sectio}
\noindent We first prove DD--4DVAR  convergence (see Figure \ref{schema1}).

\begin{figure}
\centering
{\includegraphics[width=0.9\textwidth]{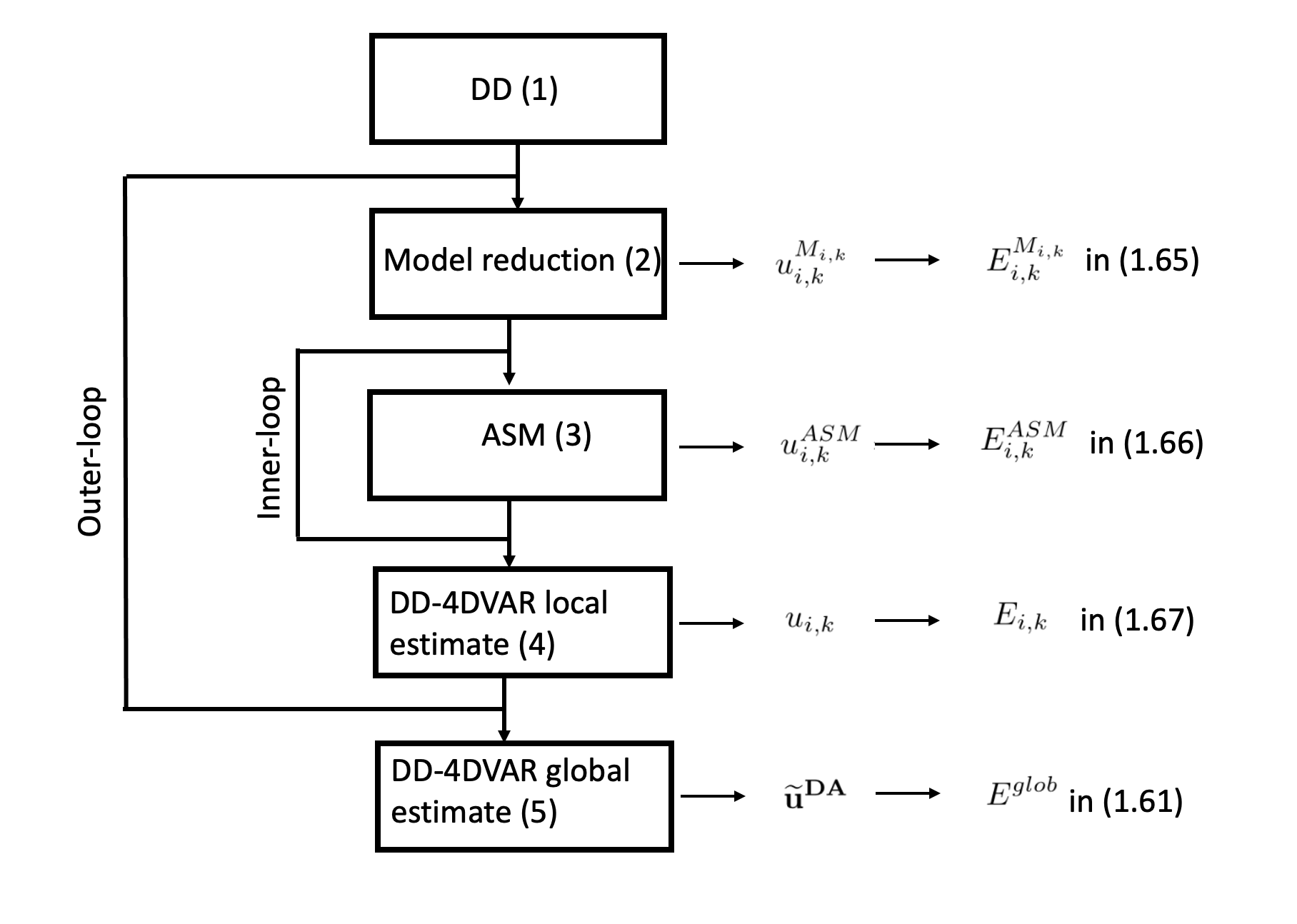}}\\
\caption{Local truncation errors related to each module of DD--4DVAR method.}
\label{schema2}
\end{figure}

\begin{theorem}\label{conv_outer}(Convergence of  outer loop) The DD-4DVAR method described in section \ref{section_DD_4DVAR} is a convergent method i.e. 
\begin{equation}\label{tesi_cov}
    lim_{n\rightarrow \infty}\|\mathbf{\Tilde{u}^{DD-DA,n+1}}-\mathbf{\Tilde{u}^{DD-DA,n}}\|=0
\end{equation}
where $\mathbf{u^{DA,n}}$ is the DD--4DVAR approximation at iteration $n$ defined in (\ref{tot_n}).
\end{theorem}
\textbf{Proof.}
From (\ref{thesis}) and (\ref{def_uDA}), it is
\begin{align}
    \begin{split}
        \mathbf{J}(\mathbf{u^{DA}})=J(\mathbf{\widetilde{u}^{DD-DA}})=\mathbf{J}(\mathbf{\widetilde{u}^{DD-DA,\bar{n}}})
    \end{split}
\end{align}
where $\bar{n}$ is the latest iteration of DD-4DVAR outer loop (see Figure \ref{schema1}). By minimum properties of $\mathbf{\widetilde{u}^{DD-DA,\bar{n}}}$, we get 
\begin{align}\label{passo2}
    \begin{split}
    \mathbf{J}(\mathbf{\widetilde{u}^{DD-DA,\bar{n}}})\le   \mathbf{J}(\mathbf{\widetilde{u}^{DD-DA,\bar{n}-1}}).
    \end{split}
\end{align}
From (\ref{passo2}), by recurring we obtain 
\begin{align*}
    \begin{split}
    0 \le \mathbf{J}(\mathbf{\widetilde{u}^{DD-DA,\bar{n}}})\le   \mathbf{J}(\mathbf{\widetilde{u}^{DD-DA,\bar{n-1}}})\le \ldots \le \mathbf{J}(\mathbf{\widetilde{u}^{DD-DA,{1}}})\le \mathbf{J}(\mathbf{\widetilde{u}^{DD-DA,{0}}})
    \end{split}
\end{align*}
then $\{\mathbf{J}(\mathbf{\widetilde{u}^{DD-DA,{n}}})\}_{\mathbf{n}\in \mathbb{N}}$ is monotonically decreasing and bounded from below by 0 then it is convergent i.e.
\begin{equation}\label{conv_dimo2}
    lim_{n\rightarrow \infty} \left[ \mathbf{J}(\mathbf{\widetilde{u}^{DD-DA,{n+1}}}) - \mathbf{J}(\mathbf{\widetilde{u}^{DD-DA,{n}}})\right]=0.
\end{equation}
Concerning iteration $n$ of the outer loop, we get 
\begin{align}
 \begin{split}
    \mathbf{J}(\mathbf{\widetilde{u}^{DD-DA,{n}}})=\|\mathbf{\widetilde{u}^{DD-DA,{n}}} -u^{{M}}\|_{\mathbf{B}^{-1}}^2+\|G \mathbf{\widetilde{u}^{DD-DA,{n}}}- v\|_{\mathbf{R}^{-1}}^2\ge \|\mathbf{\widetilde{u}^{DD-DA,{n}}} -u^{{M}}\|_{\mathbf{B}^{-1}}
     \end{split}  
\end{align}
from the triangle inequality we get 
\begin{align}\label{conv_dimos_2}
 \begin{split}
    \mathbf{J}(\mathbf{\widetilde{u}^{DD-DA,{n}}})\ge |\|\mathbf{\widetilde{u}^{DD-DA,{n}}}\|_{\mathbf{B}^{-1}}-\|u^{{M}}\|_{\mathbf{B}^{-1}}|
     \end{split}  
\end{align}
where $u^{{M}}$ is defined in (\ref{background}) and does not depend on $n$ then $\{\mathbf{\widetilde{u}^{DD-DA,{n}}}\}_{\mathbf{n}\in \mathbb{N}}$ is convergent respect to $\|\cdot\|_{\mathbf{B}^{-1}}$, consequently respect to $\|\cdot\|$. Then we get the thesis in (\ref{tesi_cov}).
\\
\begin{theorem}\label{conv_inner}(Convergence of  inner loop)
The ASM module  is  convergent.  
\end{theorem}
\noindent \textbf{Proof.} The convergence of ASM is proved in \cite{ASM}.
\\
\\
\noindent We also analyse the convergence of DD-4DVAR modules in terms of local truncation errors $E_{i,k}^{{M}_{i,k}}$, $E_{i,k}^{ASM}$, $E_{i,k}$ and $E_{g}$ reported in Figure \ref{schema2}. 
\noindent Similar to the differential equations case \cite{Dahlquist}, the local truncation errors  $E_{i,k}^{{M}_{i,k}}$, $E_{i,k}^{ASM}$, $E_{i,k}$ in $\Omega_i\times \Delta_k$ and $E_{g}$ in $\Omega\times \Delta$ are defined as the remainder after solutions $u^{\mathcal{M}}$ of model (\ref{modelloDA}) and $\mathbf{u^{DA}}$ of 4D-DA problem (\ref{varDA}) are substituted into the discrete models. To this aim, we give the following definitions.

\begin{definition}
We define 
\begin{equation}\label{sol_model}
    u^{\mathcal{M}\rightarrow M}:=M\cdot \{u^{\mathcal{M}}(x_i,t_k)\}_{\tilde{i}\in I; \ \tilde{k} \in K}
\end{equation}
approximation of $u^{\mathcal{M}}$ in (\ref{modelloDA}) in $\Omega\times \Delta$ obtained by replacing $u^{\mathcal{M}}$ evaluated in $\Omega_I\times\Delta_K$ defined in (\ref{discreto_dominio}), into  $M$ defined in (\ref{model_disc_tot}).
\end{definition}

\begin{definition}(Local truncation errors in $\Omega_i\times \Delta_k$)\label{definizione_errori}
Let us consider
\begin{itemize}
  \item $(\Omega_{1},\Omega_{2},\ldots,\Omega_{N_{sub}})$ and $(\Delta_{1},\Delta_{2},\ldots,\Delta_{N_{t}})$: be  decomposition of  $\Omega$ and $\Delta$ (DD in Figure \ref{schema2});
    \item$u^{\mathcal{M}\rightarrow M}/({\Omega_i}\times \Delta_k)$: be restriction  of approximation of $u^{\mathcal{M}}$ defined in (\ref{sol_model}) to $\Omega_i\times \Delta_k$;
    \item $u_{i,k}^{{M}_{i,k},\bar{n}}$: be solution of $P_{i,k}^{{M}_{i,k},\bar{n}}$ in $\Omega_i\times \Delta_k$ defined in (\ref{co}) at iteration $\bar{n}$ (Model Reduction in Figure \ref{schema2});
     \item $\mathbf{{u}^{DA}}$: be solution of 4DVAR DA problem in (\ref{varDA});
     \item $\mathbf{\widetilde{u}^{DD-DA}}$: be DD-4DVAR solution in $\Omega\times \Delta$ of $\mathbf{{u}^{DA}}$ in (\ref{varDA});
    \item $\mathbf{u^{DA}}/({\Omega}_{i}\times \Delta_{k})$: be the restriction to ${\Omega_{i}}\times \Delta_{k}$ of $\mathbf{{u}^{DA}}$;
    \item $u_{i,k}^{ASM,\bar{n}}$: be solution of  $P_{i,k}^{\bar{n}}$ on $\Omega_i\times \Delta_k$ defined in (\ref{co2}) (ASM in Figure \ref{schema2});
    \item $u_{i,k}^{\bar{n}}$: be local DD-4DVAR solution in $\Omega_i\times \Delta_k$  defined in (\ref{solpara}) (DD--4DVAR local solution).
\end{itemize}
Then, $\forall i=1,\ldots,N_{sub}$ and  $k=1,\ldots,N_{t}$, we define 
\begin{equation}\label{local_model_truncation_error}
    E_{i,k}^{{M}_{i,k},\bar{n}}:=\left\|u^{\mathcal{M}\rightarrow M}/({\Omega_i}\times \Delta_k)-u_{i,k}^{{M}_{i,k},\bar{n}}\right\|
\end{equation}
as local truncation error of numerical scheme restricted to $\Omega_i\times \Delta_k$ at iteration $\bar{n}$;
\begin{equation}\label{ASM_local_truncation_error}
    E_{i,k}^{ASM,\bar{n}}:=\left\|\mathbf{u^{DA}}/({\Omega}_{i}\times \Delta_{k})-u_{i,k}^{ASM,\bar{n}}\right\|
\end{equation}
as local truncation error of ASM restricted to $\Omega_i\times \Delta_k$ at iteration $\bar{n}$;
\begin{equation}\label{local_truncation_error}
    E_{i,k}^{\bar{n}}:=\left\|\mathbf{u^{DA}}/({\Omega}_{i}\times \Delta_{k})-u_{i,k}^{\bar{n}}\right\|
\end{equation}
as local truncation error of DD--4DVAR method restricted to  $\Omega_i\times \Delta_k$ at iteration $\bar{n}$;
\begin{equation}\label{tot_truncation_error}
    E_{g}:=\left\|\mathbf{u^{DA}}-\mathbf{\widetilde{u}^{DD-DA}}\right\|
\end{equation}
as local truncation error of DD--4DVAR method in $\Omega\times \Delta$.

\end{definition}

\noindent DD-4DVAR method needs few iterations of DD--4DVAR outer loop on $n$ to update DD-4DVAR approximation in (\ref{solpara}). Consequently, we neglect the dependency on outer loop iteration $\bar{n}$ of $u_{i,k}^{M_{i,k},\bar{n}}$, $u_{i,k}^{ASM,\bar{n}}$ and $u_{i,k}^{n}$ defined respectively in (\ref{co}), (\ref{co2}) and (\ref{solpara}) and $E_{i,k}^{{M}_{i,k},\bar{n}}$, $E_{i,k}^{ASM,\bar{n}}$ and $E_{i,k}^{\bar{n}}$ in $\Omega_i\times \Delta_k$  defined respectively in (\ref{local_model_truncation_error}), (\ref{ASM_local_truncation_error}) and (\ref{local_truncation_error}).\\

\noindent For $i=1,\ldots,N_{sub}$ and $k=1,\ldots,N_t$, we pose
\begin{equation}\label{zero0fun_pre}
    \mathbf{\widetilde{u}_{i,k}}:={{u}}_{i,k}^{\bar{n}}.
\end{equation}
From (\ref{zero0fun_pre}), the DD--4DVAR approximation in $\Omega\times \Delta$ defined in (\ref{def_uDA}) becomes
\begin{equation}\label{zero0fun}
\mathbf{\widetilde{u}^{DD-DA}}=\sum_{i=1}^{N_{sub}}\sum_{k=1}^{N_{t}} \mathbf{\widetilde{u}_{i,k}}^{EO}.
\end{equation}
\noindent For $i=1,\ldots,N_{sub}$ and $k=1,\ldots,N_t$, we pose
\begin{equation}\label{def_ASM_estimate}
\begin{split}
u_{i,k}^{ASM}&=u_{i,k}^{ASM,\bar{n}};\\ 
u_{i,k}^{{M}_{i,k}}&=u_{i,k}^{{M}_{i,k},\bar{n}};\\
u_{i,k}&=u_{i,k}^{\bar{n}}.
\end{split}
\end{equation}
Consequently, from  (\ref{def_ASM_estimate})  we pose 
\begin{equation}\label{model_error_no_n}
    E_{i,k}^{{M}_{i,k}}=E_{i,k}^{{M}_{i,k},\bar{n}},
\end{equation}
as local model truncation error on $\Omega_i\times \Delta_k$,  
\begin{equation}\label{ASM_error_no_n}
    E_{i,k}^{ASM}=E_{i,k}^{ASM,\bar{n}},
\end{equation}
as local ASM truncation error on $\Omega_i\times \Delta_k$ in (\ref{ASM_local_truncation_error}),
\begin{equation}\label{error_no_n}
    E_{i,k}=E_{i,k}^{\bar{n}}.
\end{equation}
as local truncation error $\Omega_i\times \Delta_k$ in (\ref{local_truncation_error}).


\noindent We introduce the definition of  consistency of DD--4DVAR method.
We pose $\|\cdot\|=\|\cdot \|_{2}$\footnote{We refer to $\|z_{i,k} \|_{2} =\|\{z_{i,k}(\bar{i},\bar{k})\}_{\bar{i}\in {I}_i,\bar{k}\in K_k}\|_2$ where $z\in R^{N_p\times N} $ and $z_{i,k}:=z/(\Omega_i\times \Delta_k)$ and ${I}_i$ and $K_k$ are defined in (\ref{insieme_I}) and (\ref{insieme_I}), respectively.}.
\begin{definition}\label{def_con}(Consistency of DD--4DVAR method)
Let $E_g$ be local truncation error of  DD--4DVAR in $\Omega\times \Delta$ defined in (\ref{tot_truncation_error}). The DD--4DVAR method is said to be consistent, if $E_g\rightarrow 0$  as $\Delta x,\ \Delta t \rightarrow 0$, where
\begin{itemize} 
  \item $\Delta x:=max_{i=1,\ldots,N_{sub}} (\Delta x)_i$, where $\{(\Delta x)_i\}_{i=1,\ldots,N_{sub}}$ are spatial step sizes of $M_{i,k}$ defined in (\ref{local_matrix});
      \item $\Delta t:=max_{k=1,\ldots,N_{t}}(\Delta t)_k$, where  $\{(\Delta t)_k\}_{k=1,\ldots,N_{t}}$ are temporal step sizes of $M_{i,k}$ defined in (\ref{local_matrix}).
\end{itemize} 
\end{definition}

\begin{figure}
\centering
{\includegraphics[width=1.05\textwidth]{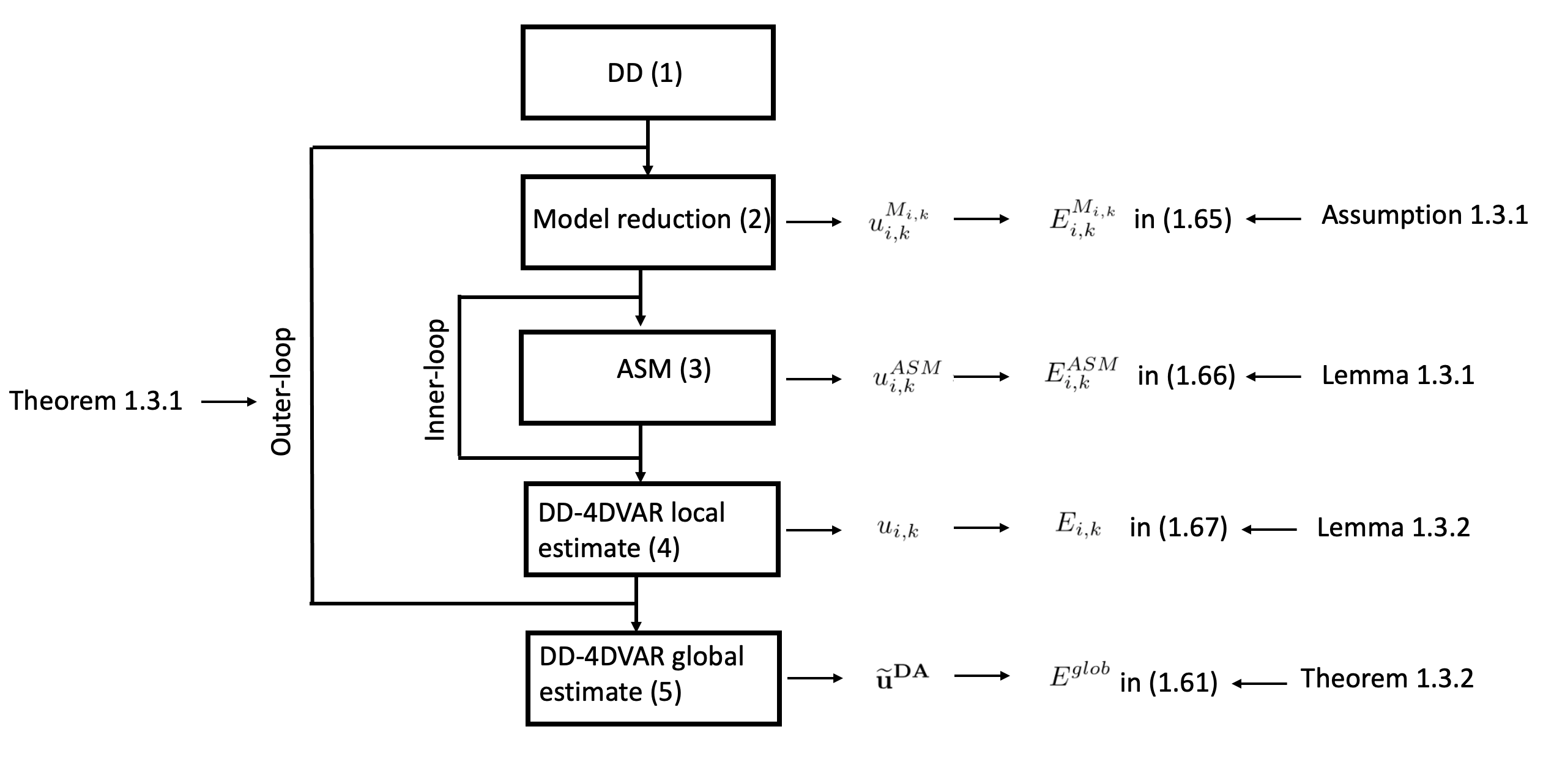}}\\
\caption{Assumption, Lemmas, Theorems related to each module of DD-4DVAR method.}
\label{schema3}
\end{figure}

 \noindent In order to prove the consistency of DD--4DVAR method, 
 we perform the analysis of local truncation errors $E_{i,k}^{{M}_{i,k}}$, $E_{i,k}^{ASM}$, $E_{i,k}$ and $E_{g}$ defined respectively  in (\ref{model_error_no_n}), (\ref{ASM_error_no_n}), (\ref{error_no_n}) and (\ref{tot_truncation_error}) related to Module Reduction, ASM, DD--4DVAR local solution and DD--4DVAR solution modules as described in Figure \ref{schema3}.\\
\begin{assumption}(Local truncation error of Model Reduction in $\Omega_i\times\Delta_k$)
Let 
\begin{equation}\label{model_truncation}
E_{i,k}^{{M}_{i,k}}=\mathcal{O}((\Delta x)_i^p+(\Delta t)_k^q), \quad \forall (i,k)\in\{1,\ldots,N_{sub}\}\times\{1,\ldots,N_t\}
\end{equation}
be local truncation error defined in (\ref{model_error_no_n}) where $(\Delta x)_i$ and $(\Delta t)_k$ are spatial and temporal step sizes of ${M}_{i,k}$ defined in (\ref{local_matrix}) and $p$ and $q$ are the order of convergence in space and in time. 
\end{assumption}
In experimental results (see section \ref{validation_section_2}), in order to discretize the Shallow Water Equations (SWEs) model we consider Lax--Wendroff scheme \cite{LeVeque}. Hence, in that case, $p=q=2$.

\begin{lemma}\label{lemma_ASM}(Local truncation error of ASM in $\Omega_i\times\Delta_k$) Let us consider 
\begin{itemize}
\item $(\Omega_{1},\Omega_{2},\ldots,\Omega_{N_{sub}})$ and $(\Delta_{1},\Delta_{2},\ldots,\Delta_{N_{t}})$: be  decomposition of  $\Omega$ and $\Delta$ (DD in Figure \ref{schema3});
\item $E_{i,k}^{ASM}$: be local truncation error in $\Omega_i\times \Delta_k$ of ASM defined in (\ref{ASM_error_no_n}); 
\item $u_{i,k}^{ASM}$: be solution of $P_{i,k}^{\bar{n}}$ on $\Omega_i\times \Delta_k$ defined in (\ref{def_ASM_estimate}) (ASM in Figure \ref{schema3});
    \item $\sigma_0^2$: be observational error variance;
    \item $\mathbf{B}_i=\mathbf{V}_i\mathbf{V}_i^T$: be restriction of covariance matrices of the error on background to $\Omega_{i}$;
    \item $G_{i,k}$: be restriction to $\Omega_i\times \Delta_k$ of matrix $G$ defined in (\ref{mat_G});
    \item $ad_i$: be number of subdomains adjacent to $\Omega_i$, defined in (\ref{ni});
    \item $\mathbf{B}_{ij}=\mathbf{V}_{ij}\mathbf{V}_{ij}^T$: be restriction to $\Omega_{ij}$ of variance  matrix of the error on background;
    \item ${M}_{i,k}$: be matrix defined in (\ref{matriceM}), restricted to $\Omega_i\times \Delta_k$;
    \item $\mu(\mathbf{V}_i)$, $\mu(G_{i,k})$, $\mu(M_{i,k})$ and $\mu(\mathbf{V}_{ij})$: be condition number of matrices $\mathbf{V}_i, G_{i,k}$, $M_{i,k}$ and $\mathbf{V}_{ij}$, respectively.
\end{itemize}
Then, $\forall i=1,\ldots,N_{sub}; k=1,\ldots,N_t$, it holds that:
\begin{equation}\label{tesi_conv}
   E_{i,k}^{ASM}\le \mu_{i,k}^{DD-DA}\cdot E_{i,1}^{ASM}
\end{equation}
where
\begin{align}\label{def_parametri}
\begin{split}
\mu_{i,k}^{DD-DA}&:=\left[1+\frac{1}{\sigma_0^2}\mu^2(\mathbf{V}_i)\mu^2(G_{i,k})+ad_i\cdot \mu^2(\mathbf{V}_{ij})\right]\mu({M}_{i,k}).
\end{split}
\end{align}
\end{lemma}
\textbf{Proof:}
By applying error approximation in \cite{JCP}, we have
\begin{equation}\label{1_t}
\|\mathbf{u}^{DA}/(\Omega_i\times \Delta_k)-u_{i,k}^{ASM}\|\le \mu(\mathbf{J}_{i,k})\mu({M}_{i,k})\cdot \|\mathbf{u^{DA}}/(\Omega_i\times \Delta_1)-u_{i,1}^{ASM}\|
\end{equation}
where $\|\mathbf{u^{DA}}/(\Omega_i\times \Delta_1)-u_{i,1}^{ASM}\|$ is error on $\Omega_i\times \Delta_1$. As proved in \cite{JCP}, it is
\begin{equation}
    \mu(\mathbf{J}_{i,k})=\mu(A_{i,k})
\end{equation}
where $\mathbf{J}_{i,k}$ and $A_{i,k}$  are respectively defined in (\ref{local_operator}) and (\ref{Aik}) and form the triangle inequality, it is 
\begin{align}\label{2_t}
    \mu(A_{i,k})&\le 1+\frac{1}{\sigma_0^2}\mu^2(\mathbf{V}_i)\mu^2(G_{i,k})+ad_i\cdot \mu(\mathbf{B}_{ij})\\
&\le  1+\frac{1}{\sigma_0^2}\mu^2(\mathbf{V}_i)\mu^2(G_{i,k})+ad_i\cdot \mu^2(\mathbf{V}_{ij}),
\end{align}
where $\sigma_0^2$ is observational error variance, $ad_i$ is number of adjacent subdomains to $\Omega_i$ defined in (\ref{ni}), $\mu(\mathbf{V}_i)$, $\mu(G_{i,k})$, $\mu(M_{i,k})$ and $\mu(\mathbf{V}_{ij})$ are condition numbers of matrices $\mathbf{V}_i, G_{i,k}$, $M_{i,k}$ and $\mathbf{V}_{ij}$, respectively.
From (\ref{1_t}) and (\ref{2_t}), the (\ref{tesi_conv}) follows. 

\noindent In the following we approximation the local truncation error of DD-4DVAR method in $\Omega_i\times \Delta_k$ defined in (\ref{error_no_n}). 
\begin{lemma}(Local truncation error of DD-4DVAR method in $\Omega_i\times \Delta_k$)\label{convergenza}
Let us consider
\begin{itemize}
    \item $(\Omega_{1},\Omega_{2},\ldots,\Omega_{N_{sub}})$ and $(\Delta_{1},\Delta_{2},\ldots,\Delta_{N_{t}})$: be decomposition of  $\Omega$ and $\Delta$ (DD in Figure \ref{schema3});
  \item $\mathbf{u^{DA}}/(\Omega_{i}\times \Delta_{k})$: be the restriction to $\Omega_{i}\times \Delta_{k}$ of 4DVAR DA problem's solution in (\ref{varDA});
   \item $u_{i,k}$: be DD-4DVAR local solution on $\Omega_i\times \Delta_k$ at iteration ${n}$ defined in (\ref{solpara}) (DD-4DVAR local solution in Figure \ref{schema3}); 
    \item$u^{\mathcal{M}\rightarrow M}/({\Omega_i}\times \Delta_k)$: be restriction of  $u^{\mathcal{M}}$ defined in (\ref{sol_model}) to $\Omega_i\times \Delta_k$;
\item $E_{i,k}^{{M}_{i,k}}$: be local truncation error in $\Omega_i\times \Delta_k$ defined in (\ref{model_error_no_n});
  \item $E_{i,1}^{ASM}$: be local truncation error in $\Omega_i\times \Delta_1$ defined in (\ref{ASM_error_no_n});
\end{itemize}
Then $\forall i=1,\ldots,N_{sub};k=1,\ldots,N_{t}$, it holds that:
\begin{align}\label{tesi1}
\begin{split}
E_{i,k} \le \mu_{i,k}^{DD-DA}\cdot E_{i,1}^{ASM}+2\cdot  E_{i,k}^{{M}_{i,k}}
\end{split}
\end{align}
where $\mu_{i,k}^{DD-DA}$ is defined in (\ref{def_parametri}).
\end{lemma}
\noindent \textbf{Proof}: From (\ref{def_ASM_estimate}) and (\ref{solpara}), $E_{i,k}$ defined in (\ref{error_no_n}) can be rewritten as follows
\begin{equation*}
\begin{array}{ll}
E_{i,k}&:=\|\mathbf{u^{DA}}/(\Omega_{i}\times \Delta_{k})-u_{i,k}\|\\&=\|\mathbf{u^{DA}}/(\Omega_{i}\times \Delta_{k})-u_{i,k}^{{M}_{i,k},\bar{n}}-(u_{i,k}^{ASM,\bar{n}-1}-u_{i,k}^{{M}_{i,k},\bar{n}-1})\|
\end{array}
\end{equation*}
from the triangle inequality
\begin{equation}\label{pri}
\begin{array}{ll}
E_{i,k}&\le \|\mathbf{u^{DA}}/(\Omega_{i}\times \Delta_{k})-u_{i,k}^{ASM,\bar{n}}\|+\|u_{i,k}^{{M}_{i,k},\bar{n}-1}-u_{i,k}^{{M}_{i,k},\bar{n}}\|
\end{array}
\end{equation}
as consequence of Lemma \ref{lemma_ASM} and (\ref{ASM_error_no_n}) we have 
\begin{equation}\label{dim}
 E_{i,k}\le  \mu_{i,k}^{DD-DA}\cdot E_{i,1}^{ASM}+\|u_{i,k}^{{M}_{i,k},\bar{n}-1}-u_{i,k}^{{M}_{i,k},\bar{n}}\|
\end{equation}
where $E_{i,1}^{ASM}$ is defined in (\ref{ASM_error_no_n}) and $\mu_{i,k}^{DD-DA}$ is defined in (\ref{def_parametri}).\\
In particular, by adding and subtracting $u^{\mathcal{M}\rightarrow M}/({\Omega_i}\times \Delta_k)$ in $\|u_{i,k}^{{M}_{i,k},\bar{n}-1}-u_{i,k}^{{M}_{i,k},\bar{n}}\|$ we get:
\begin{align*}
\begin{split}
    \|u_{i,k}^{{M}_{i,k},\bar{n}-1}-u_{i,k}^{{M}_{i,k},\bar{n}}\|&= \|(u_{i,k}^{{M}_{i,k},\bar{n}-1}-u^{\mathcal{M}\rightarrow M}/({\Omega_i}\times \Delta_k))\\&+(u^{\mathcal{M}\rightarrow M}/({\Omega_i}\times \Delta_k)-u_{i,k}^{{M}_{i,k},\bar{n}})\|
    \end{split}
\end{align*}
and from the triangle inequality
\begin{align}\label{dim1}
\begin{split}
    \|u_{i,k}^{{M}_{i,k},\bar{n}-1}-u_{i,k}^{{M}_{i,k},\bar{n}}\| & \le \|u_{i,k}^{{M}_{i,k},\bar{n}-1}-u^{\mathcal{M}\rightarrow M}/({\Omega_i}\times \Delta_k)\|\\&+\|u^{\mathcal{M}\rightarrow M}/({\Omega_i}\times \Delta_k)-u_{i,k}^{{M}_{i,k},\bar{n}}\|\\
    &=E_{i,k}^{{M}_{i,k},\bar{n}}+E_{i,k}^{{M}_{i,k},\bar{n}-1}
    \end{split}
\end{align}
From  Theorem \ref{conv_outer}, we get that $\{\mathbf{\widetilde{u}^{DD-DA,n}}\}_{n\in \mathbb{N}}$ is a convergent sequence, then it is a Cauchy sequence. From (\ref{tot_n}), we get that $\{u_{i,k}^{M_{i,k},n}\}_{n\in \mathbb{N}}$ is also a Cauchy sequence, i.e.
\begin{equation}\label{cauchy_conv}
    \forall \  \epsilon>0 \ \exists \ \mathbf{N}>0 \ : \ \ \|u_{i,k}^{M_{i,k},{n}}- u_{i,k}^{M_{i,k},m}\|\le \epsilon \quad \forall n,m>\mathbf{N}.
\end{equation}
In particular, (\ref{cauchy_conv}) is true for $n=\bar{n}$ and $m=\bar{n}-1$, assuming that $\bar{n}$ is large enough. Consequently, we can neglect the dependency on outer loop in $E_{i,k}^{{M}_{i,k},\bar{n}}$ and $E_{i,k}^{{M}_{i,k},\bar{n}-1}$ in (\ref{dim1}) i.e.
\begin{align}\label{dim2}
    \|u_{i,k}^{{M}_{i,k},\bar{n}-1}-u_{i,k}^{{M}_{i,k},\bar{n}}\|\le
    2 \cdot E_{i,k}^{{M}_{i,k}}
\end{align}
where $E_{i,k}^{{M}_{i,k}}$ is defined in (\ref{model_error_no_n}). From (\ref{dim}), (\ref{dim1}) and (\ref{dim2}) we get the thesis in (\ref{tesi1}).
\\
\normalsize

\begin{lemma}\label{lemma_error}
If $e_0$, the  error  on initial condition of ${M}$ in (\ref{model_disc_tot}) is equal to zero i.e. 
\begin{equation}\label{ipotesi_lemma}
e_0=0    
\end{equation}
then
\begin{equation}\label{local_trunc_error}
    E_{i,k}\le c\cdot ((\Delta x)_i^p+\cdot (\Delta t)_k^q)
\end{equation}
where $E_{i,k}$ is the local truncation error defined in (\ref{error_no_n}) and $c$ is positive constant independent on DD.
\end{lemma}
\noindent 
\textbf{Proof.} By applying the Lemma \ref{ASM_error_no_n} to local ASM truncation error $E_{i,1}^{ASM}$ on $\Omega_i\times \Delta_1$, we get 
\begin{equation}\label{parte_1_ASM}
  E_{i,1}^{ASM}\le  F_{i,1} \cdot e_0/\Omega_i. 
\end{equation}
where $F_{i,1}$ is defined in (\ref{def_parametri}) and $ e_0/\Omega_i$ is the restriction of $e_0$ to $\Omega_i$.
From the assumptions it is $e_0/\Omega_i=0$ and by replacing it in (\ref{parte_1_ASM}) we have
\begin{equation}\label{parte_1_b_ASM}
 E_{i,1}^{ASM}=0   
\end{equation}
consequently local ASM truncation error on $\Omega_i\times \Delta_k$  is 
\begin{equation}\label{parte_2_ASM}
 E_{i,k}^{ASM}=0.   
\end{equation}
 From (\ref{parte_2_ASM}), (\ref{model_truncation}) and (\ref{tesi1}) we get the thesis in (\ref{local_trunc_error}).
\\
\\
\noindent
Now, we are able to prove the following result on global truncation error. 
\begin{theorem}\label{dd_truncation_error}(Local truncation error in $\Omega\times \Delta$) Let assume the assumption of Lemma \ref{lemma_error} in (\ref{ipotesi_lemma}). The local truncation error in $\Omega\times \Delta$ is 
\begin{equation}\label{global_error}
    E_{g}\le c\cdot (N_{sub} \cdot N_{t}) \cdot [(\Delta x)^p+ (\Delta t)^q],
\end{equation}
 where 
 \begin{itemize}
     \item $\mathbf{u^{DA}}$: solution of 4DVAR DA problem in (\ref{varDA});
     \item $\mathbf{\widetilde{u}^{DD-DA}}$: DD--4DVAR solution in $\Omega\times \Delta$ defined in (\ref{def_uDA});
  \item $\Delta x:=max_{i=1,\ldots,N_{sub}} (\Delta x)_i$: where $\{(\Delta x)_i\}_{i=1,\ldots,N_{sub}}$ are spatial step sizes of $M_{i,k}$ defined in (\ref{local_matrix});
      \item $\Delta t:=max_{k=1,\ldots,N_{t}}(\Delta t)_k$: where  $\{(\Delta t)_k\}_{k=1,\ldots,N_{t}}$ are temporal step sizes of $M_{i,k}$ defined in (\ref{local_matrix});
      \item $c>0$: is positive constant independent of DD.
     
  \end{itemize}
\end{theorem}

\textbf{Proof} From (\ref{zero0fun}) the error $E_{g}$ defined in (\ref{tot_truncation_error}) can be rewritten as follows 
\begin{align*}
 \begin{split}
 E_{g}:=\left\|\mathbf{u^{DA}}-\mathbf{\tilde{u}^{DD-DA}}\right\|=\left\|\mathbf{u^{DA}}-\sum_{i=1}^{N_{sub}}\sum_{k=1}^{N_t}{\mathbf{\widetilde{u}_{i,k}}}\right\|
 \end{split}
\end{align*}
by applying the restriction and extension operator (Definitions \ref{def_rest} and \ref{def_ex}) to $\mathbf{u^{DA}}$, we get 
\begin{align*}
 \begin{split}
 E_{g}=\left\| \sum_{i=1}^{N_{sub}}\sum_{k=1}^{N_t}\left[ (\mathbf{u^{DA}}/(\Omega_i\times \Delta_k))^{EO}-{\mathbf{\widetilde{u}_{i,k}}}\right]\right\|
 \end{split}
\end{align*}
from triangle inequality, it is 
\begin{align}\label{dim_tot_1}
 \begin{split}
 E_{g}=\left\| \sum_{i=1}^{N_{sub}}\sum_{k=1}^{N_t}\left[ (\mathbf{u^{DA}}/(\Omega_i\times \Delta_k))^{EO}-{\mathbf{\widetilde{u}_{i,k}}}\right]\right\|\\\le \sum_{i=1}^{N_{sub}}\sum_{k=1}^{N_t}\left\|    (\mathbf{u_{i,k}^{DA}})^{EO}-{\mathbf{\widetilde{u}_{i,k}}}\right\|=\sum_{i=1}^{N_{sub}}\sum_{k=1}^{N_t}E_{i,k}
 \end{split}
\end{align}
where $E_{i,k}$ is defined in (\ref{error_no_n}).
From Lemma \ref{lemma_error} we have
\begin{align}\label{dim_error_1}
 \begin{split}
  \sum_{i=1}^{N_{sub}}\sum_{k=1}^{N_t}E_{i,k}&\le   c\cdot  \sum_{i=1}^{N_{sub}}\sum_{k=1}^{N_t}  ((\Delta x)_i^p+\cdot (\Delta t)_k^q)\\
  &=c\cdot \left[N_t\cdot \sum_{i=1}^{N_{sub}} (\Delta x)_i^p+N_{sub}\sum_{k=1}^{N_t} (\Delta t)_k^q\right]
   \end{split}  
\end{align}
consequently
\begin{align*}
 \begin{split}
E_{g}\le c\cdot \left[N_t\cdot \sum_{i=1}^{N_{sub}} (\Delta x)_i^p+N_{sub}\sum_{k=1}^{N_t} (\Delta t)_k^q\right].
 \end{split}
\end{align*}
By defining
\begin{align*}
    \begin{split}
 \Delta x&:=max_{i=1,\ldots,N_{sub}} (\Delta x)_i\\
\Delta t&:=max_{k=1,\ldots,N_{t}}(\Delta t)_k
    \end{split}
\end{align*}
we get 
\begin{align*}
 \begin{split}
E_{g}&\le c\cdot \left[N_t\cdot \sum_{i=1}^{N_{sub}} (\Delta x)_i^p+N_{sub}\sum_{k=1}^{N_t} (\Delta t)_k^q\right]\\
&\le  c\cdot \left[N_t\cdot \sum_{i=1}^{N_{sub}} (\Delta x)^p+N_{sub}\sum_{k=1}^{N_t} (\Delta t)^q\right]\\
&\le   c\cdot (N_{sub} \cdot N_{t}) \cdot [(\Delta x)^p+ (\Delta t)^q].
 \end{split}
\end{align*}
Hence, the (\ref{global_error}) is proved.\\

\begin{figure}
\centering
{\includegraphics[width=1.2\textwidth]{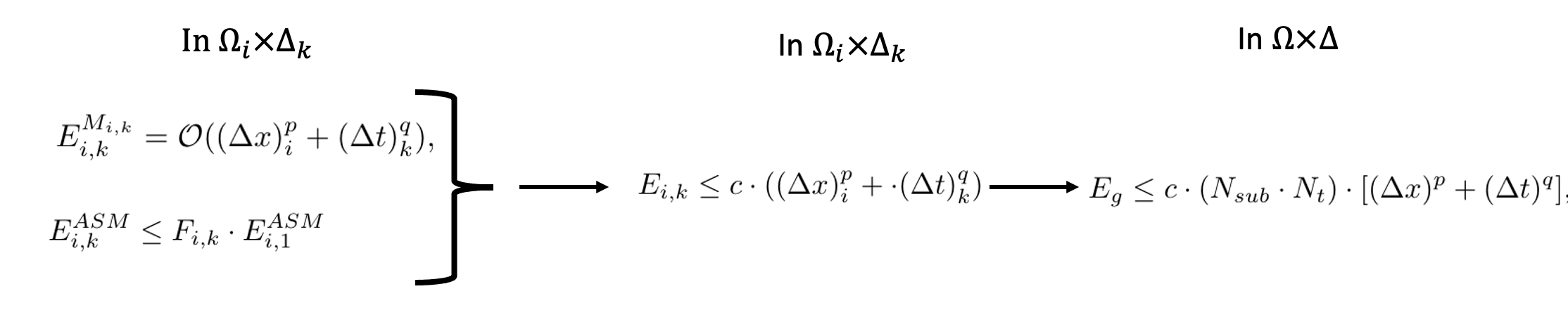}}\\
\caption{Schematic description of DD--4DVAR local truncation analysis.}
\label{schema_finale}
\end{figure}

\noindent In Figure \ref{schema_finale}, it is schematic description of results obtained from DD--4DVAR  local truncation analysis. \\

\noindent Now, we prove the stability of DD--4DVAR method by studying the propagation error from each time interval to the next, assuming that the predictive model is stable.

\begin{assumption}\label{ipotesi_model}(Stability of model discrete scheme) The discrete scheme applied to $\mathcal{M}$ in (\ref{modelloDA}) is stable scheme i.e. $\exists D>0 $ such that
\begin{equation}
 \|u^{M}-v^{M}\|\le D\cdot e_0,
\end{equation}
where
\begin{itemize}
    \item $u^{M}$: is solution of $M$ in (\ref{model_disc_tot});
    \item $v^{M}$: is solution of $\bar{M}$, where $\bar{M}$ is obtained by considering initial error $e_0$ on initial condition;
    \item $e_0$: is initial error on initial condition of $M$ in (\ref{model_disc_tot}).
\end{itemize}
\end{assumption}
\begin{definition}(Propagation error from $\Delta_{k-1}$ to $\Delta_{k}$)
Let us consider 
\begin{itemize}
    \item $(\Delta_1,\Delta_2,\ldots,\Delta_{N_t})$: be the decomposition of $\Delta$ (DD module);
     \item $\mathbf{\widetilde{u}^{DD-DA}}$: be DD--4DVAR solution in $\Omega\times \Delta$ defined in (\ref{def_uDA});
    \item ${\mathbf{\Tilde{v}^{DD-DA}}}$: be DD--4DVAR solution in $\Omega\times \Delta$ computed by adding perturbation $e_k$ to initial condition of $P_{i,k}^{M_{i,k},n}$ defined in (\ref{co}).
\end{itemize}
We define 
\begin{equation}\label{def_prop_error}
  \bar{E}_k:=\|\mathbf{\Tilde{u}^{DD-DA}}/\Delta_k-{\mathbf{\Tilde{v}^{DD-DA}}}/\Delta_k\|   
\end{equation}
the propagation error from $\Delta_{k-1}$ to $\Delta_{k}$.
\end{definition}

\noindent In the following, result of stability of DD--4DVAR method.
\begin{theorem}\label{teo_stab}(Stability of DD--4DVAR method)\label{stabilita}
If $e_0$, the error  on initial condition of ${M}$ in (\ref{model_disc_tot}), is equal to zero i.e. 
\begin{equation}\label{ipotesi_teo}
e_0=0    
\end{equation}
then, $\forall k=1,\ldots,N_t$ $\exists \ C_k>0$ such that
\begin{equation*}
    \bar{E}_{k}\le C_k \cdot \bar{e}_k
\end{equation*}
where
\begin{itemize}
\item $\bar{E}_{k}$: is propagation error from $\Delta_{k-1}$ to $\Delta_k$ defined in (\ref{def_prop_error});
    \item $C_k$: is a constant dependind on the  model reduction and on ASM modules of DD--4DVAR method;
    \item $\bar{e}_k$: perturbation on initial condition of $P_{i,k}^{M_{i,k}}$ defined in (\ref{co}).
    \end{itemize}
\end{theorem}
\noindent \textbf{Proof.} To simplify the notations in the proof, we consider $\bar{e}_k=\bar{e}, \ \forall k=1,\ldots,N_t$. From (\ref{solpara}), (\ref{def_ASM_estimate}) and triangle inequality, we get 
\begin{align}
    \begin{split}
        \bar{E}_k:=\|\mathbf{\Tilde{u}^{DD-DA}}/\Delta_k-{\mathbf{\Tilde{v}^{DD-DA}}}/\Delta_k\|&\le \|(u_{i,k}^{M_{i,k},\bar{n}})^{EO}/\Delta_k-(v_{i,k}^{M_{i,k},\bar{n}})^{EO}/\Delta_k\|\\&+\|(u_{i,k}^{M_{i,k},\bar{n}-1})^{EO}/\Delta_k-(v_{i,k}^{M_{i,k},\bar{n}-1})^{EO}/\Delta_k\|\\&+\|(u_{i,k}^{ASM,\bar{n}})^{EO}/\Delta_k-(v_{i,k}^{ASM,\bar{n}})^{EO}/\Delta_k\|. 
    \end{split}
\end{align}
From (\ref{cauchy_conv}) and (\ref{ASM_error_no_n}), we can neglect the dependency on outer loop iteration $\bar{n}$ i.e.
\begin{align}\label{proof_stab_1}
    \begin{split}
        \bar{E}_k\le& 2\cdot \|(u_{i,k}^{M_{i,k}})^{EO}/\Delta_k-(v_{i,k}^{M_{i,k}})^{EO}/\Delta_k\|+\|(u_{i,k}^{ASM})^{EO}/\Delta_k-(v_{i,k}^{ASM})^{EO}/\Delta_k\|. 
    \end{split}
\end{align}
From Assumption \ref{ipotesi_model}, we get that $\exists \bar{D}>0$ such that
\begin{align}\label{proof_stab_2}
    \begin{split}
     \|(u_{i,k}^{M_{i,k}})^{EO}/\Delta_k-(v_{i,k}^{M_{i,k}})^{EO}/\Delta_k\|\le  \bar{D}\cdot e_{0}.
    \end{split}
\end{align}
where $e_0$ is the error on initial condition of ${M}$ in (\ref{model_disc_tot}).
By adding and subtracting $\mathbf{u^{DA}}/ \Delta_k$ to $\left[(u_{i,k}^{ASM})^{EO}/\Delta_k-(v_{i,k}^{ASM})^{EO}/\Delta_k\right]$ and from triangle inequality, it is 
\begin{align}\label{prov_sens}
    \begin{split}
    \|(u_{i,k}^{ASM})^{EO}/\Delta_k-(v_{i,k}^{ASM})^{EO}/\Delta_k\|\le    \|\mathbf{u^{DA}}/ \Delta_k-(u_{i,k}^{ASM})^{EO}/\Delta_k\|\\+\|\mathbf{u^{DA}}/ \Delta_k-(v_{i,k}^{ASM})^{EO}/\Delta_k\|.
    \end{split}
\end{align}
From (\ref{prov_sens}) and (\ref{tesi_conv}), we get
\begin{align}\label{proof_stab_3}
    \begin{split}
  \|(u_{i,k}^{ASM})^{EO}/\Delta_k-(v_{i,k}^{ASM})^{EO}/\Delta_k\|\le \mu_k^{DD-DA}\cdot E_1^{ASM}+\bar{\mu}_k^{DD-DA}\cdot \bar{E}_1^{ASM} 
    \end{split}
\end{align}
where 
\begin{align}\label{E_1}
    \begin{split}
     E_1^{ASM}&= \|\mathbf{u^{DA}}/ \Delta_1-(u_{i,1}^{ASM})^{EO}/\Delta_1\|   \\
     \bar{E}_1^{ASM}&=\|\mathbf{u^{DA}}/ \Delta_1-(v_{i,1}^{ASM})^{EO}/\Delta_1\|.
    \end{split}
\end{align}
and 
\begin{align}\label{mu}
\begin{split}
\mu_k^{DD-DA}&:=\left[1+\frac{1}{\sigma_0^2}\mu^2(\mathbf{V})\mu^2(G/\Delta_k)\right]\mu({M}/\Delta_k)\\
\bar{\mu}_k^{DD-DA}&:=\left[1+\frac{1}{\sigma_0^2}\mu^2(\mathbf{V})\mu^2(G/\Delta_k)\right]\mu(\bar{M}/\Delta_k)\\
\end{split}
\end{align}
with $\sigma_0$ observational error variance, $\mathbf{B}=\mathbf{V}\mathbf{V^T}$ covariance matrix of the error on background to $\Omega$, $G$ matrix defined in (\ref{mat_G}), $M$ defines in (\ref{model_disc_tot})  and $\bar{M}$ is discrete model obtained by considering initial error $e_0$ on initial condition of $M$.
By applying (\ref{tesi_conv}) to $E_1^{ASM}$ and $\bar{E}_1^{ASM}$ in (\ref{E_1}), we get
\begin{align}\label{prov_sens_2}
    \begin{split}
  \|(u_{i,k}^{ASM})^{EO}/\Delta_k-(v_{i,k}^{ASM})^{EO}/\Delta_k\|\le \mu_k^{DD-DA}\cdot \mu_1^{DD-DA} {e}_0+\bar{\mu}_k^{DD-DA}\cdot \bar{\mu}_k^{DD-DA} \bar{e}. 
    \end{split}
\end{align}
From (\ref{proof_stab_1}), (\ref{proof_stab_2}) and (\ref{proof_stab_3}), it is 
\begin{align*}
    \begin{split}
        \bar{E}_{k}\le 2\cdot \bar{D}\cdot e_0+\mu_k^{DD-DA}\cdot \mu_1^{DD-DA} {e}_0+\bar{\mu}_k^{DD-DA}\cdot \bar{\mu}_1^{DD-DA} \bar{e}  
    \end{split}
\end{align*}
and from hypothesis in (\ref{ipotesi_teo}), we get
\begin{align}\label{tesi_stab}
    \begin{split}
        \bar{E}_{k}\le \bar{\mu}_k^{DD-DA}\cdot \bar{\mu}_1^{DD-DA} \bar{e}  
    \end{split}
\end{align}
Consequently, for $k=1,\ldots,N_t$ we get that $\exists \ C_k>0$ such that
\begin{align}
    \begin{split}
        \bar{E}_k\le C_k\cdot \bar{e}
    \end{split}
\end{align}
where \begin{equation}\label{costante_C}
   C_k :=\bar{\mu}_k^{DD-DA}\cdot \bar{\mu}_1^{DD-DA}.
\end{equation}
The thesis is proved.
\\
\\
From Theorem \ref{teo_stab}, we get the stability of DD--4DVAR method.\\
\\
\noindent \textbf{Remark (Condition number of DD--4DVAR method)} \textit{We note that in Theorem \ref{teo_stab}, to prove the stability of DD--4DVAR method, we study the propagation error according to  forward error analysis. In particular, we get that the constant $C_k,\ \forall k=1,\ldots,N_t$ in (\ref{costante_C}) depends on the quantity $\bar{\mu}_k^{DD-DA}$ defined in (\ref{mu}), which is in turn depends on condition numbers of $M$ in (\ref{model_disc_tot}), $G$ in (\ref{mat_G}) and $\mathbf{B}=\mathbf{V}\mathbf{V^T}$ covariance matrix of the error
on background in $\Omega$. As a consequence of the forward error analysis, we may say that the quantity $\bar{\mu}_k^{DD-DA}$ can be regarded as the  condition  number of DD--4DVAR method.}
\\

\section{Validation analysis}\label{validation_section_2}

\noindent Validation is performed on the high performance hybrid computing architecture of the SCoPE (Sistema Cooperativo Per Elaborazioni scientifiche multidiscipliari) data center, located in the University of Naples Federico II. Specifically, the HPC architecture is composed by 8 nodes which consist of distributed memory DELL M600 blades. The blades are connected by a 10 Gigabit Ethernet technology and each of them is composed of 2 Intel Xeon@2.33GHz quadcore processors sharing the same local 16 GB RAM memory for a number of 8 cores per blade and of 64 total cores. Experimental results allows to deduce  that the experimental order of DD--4DVAR consistency corresponds to the theoretical one obtained in Theorem \ref{dd_truncation_error} and DD--4DVAR is well-conditioned. 
\noindent We consider the following experimental set up.\\

\noindent \textbf{4DVAR DA set up}.
\begin{itemize}
\item $\Omega=(0,1)\subset \mathbb{R}$:  spatial domain;
\item $\Delta=[0,1.5]\subset \mathbb{R}$:  time interval;
\item $N_{p}$: numbers of inner nodes of $\Omega$ defined in (\ref{setI});
\item  $N$: numbers of instants of time in $\Delta$;
\item $n_{obs}=64$:  number of observations considered at each step $l=0,1,\ldots,N$;
\item $y\in \mathbb{R}^{N\cdot n_{obs}}$: observations vector at each step $l=0,1,\ldots,N$;
\item $H_{l}\in \mathbb{R}^{n_{obs}\times N_{p}}$: piecewise linear interpolation operator whose coefficients are computed using the
nodes of $\Omega$ nearest to the observation values;
\item $G\in \mathbb{R}^{N\cdot n_{obs}\times N_{p}}$: obtained as in (\ref{mat_G}) from the matrix $H_{l}$, $l=0,1,\ldots,N$; 
\item $\sigma_{m}^{2}=0.5$, $\sigma_{0}^{2}=0.5$: model and observational error variances;
\item $\textbf{B}\equiv B_{l}=\sigma_{m}^{2}\cdot C$: covariance matrix of the error of the model at each step $l=0,1,\ldots,N$, where $C\in \mathbb{R}^{N_{p}\times N_{p}}$ denotes the Gaussian correlation structure of the model errors in (\ref{matC});
\item $ {R}_{l}=\sigma_{0}^{2}\cdot I_{n_{obs},n_{obs}}\in \mathbb{R}^{n_{obs}\times n_{obs}}$: covariance matrix of the errors of the observations at each step $l=0,1,\ldots,N-1$.
\item $\textbf{R}\in \mathbb{R}^{N\cdot n_{obs}\times N\cdot n_{obs}}$: a diagonal matrix obtained from the matrices $R_{l}$, $l=0,1,\ldots,N-1$.
\end{itemize}

\begin{figure}
\centering
{\includegraphics[width=1.\textwidth]{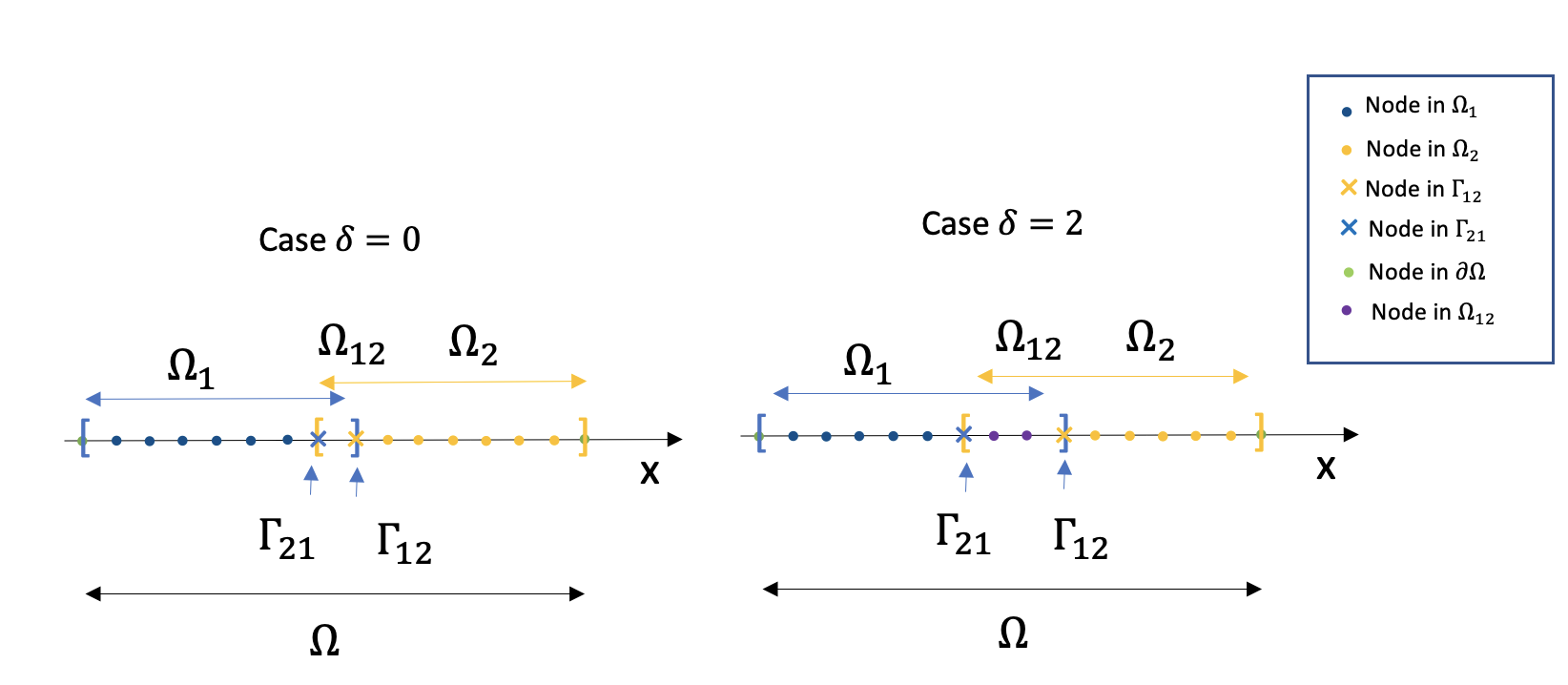}}
{\caption{Decomposition of spatial domain $\Omega\subset \mathbb{R}$ in two subdomains $\{\Omega_i\}_{i=1,2}$ by identifying overlap region $\Omega_{12}$ defined in (\ref{overlap_region}) and interfaces $\Gamma_{12}$ and $\Gamma_{21}$ defined in (\ref{interfacce}). On the left case $\delta=0$ i.e. no inner nodes in $\Omega_{12}$, on the right case $\delta=2$ i.e. two inner nodes in overlap region $\Omega_{12}$.}
\label{fig_DD_1D}}
\end{figure}

\noindent \textbf{DD-4DVAR set up}: we consider the following set up:
\begin{itemize}
\item $p$:  number of cores;
\item $N_{sub}\equiv p$: number of spatial subdomains;
\item $N_{t}$: number of time intervals;
\item $n_{1}=n_{nsub}=1$ and $n_{2}=n_3=\ldots=n_{nsub-1}=2$: number of subdomains adjacent to $\Omega_{1}, \ \Omega_{nsub}$ and $\Omega_{2},\ \Omega_{3},\ldots,\Omega_{nsub-1}$, respectively; 
\item $\delta$: number of inner nodes of overlap regions defined in (\ref{delta});
\item $N_{loc}$: inner nodes of subdomains defined in (\ref{nloc1});
\item $\Delta x$ and $\Delta t$: spatial and temporal step sizes of $M_{i,k}$ defined in (\ref{local_matrix});
\item $C:=\{c_{i,j}\}_{i,j=1,\ldots,N_{p}}\in \mathbb{R}^{N_{p}\times N_{p}}$: the Gaussian correlation structure of the model error where
 \begin{equation}\label{matC}
c_{i,j}=\rho^{|i-j|^{2}}, \quad \rho=exp\left(\frac{-\Delta x^{2}
}{2}\right), \quad |i-j|<N_{p}/2 \quad  \begin{array}{ll}\textrm{for}\ i,j=1,\ldots,{N_{p}}  \end{array}.
\end{equation}
\end{itemize}
\noindent 
We introduce 
\begin{equation}\label{error_exp}
e^p=e^p(\delta,\Delta x,\Delta t, N_{sub},N_t):=\|\mathbf{u^{DA}}-\mathbf{\Tilde{u}^{DD-DA}}\|_{2},
\end{equation}
where  $\mathbf{u^{DA}}$ denotes the minimum of the
4DVAR (global) functional $\mathbf{J}$ in (\ref{funzionale}) while $\mathbf{\Tilde{u}^{DD-DA}}$ is obtained by gathering minimum of the local 4DVar functionals $\mathbf{J}_{i,k}$ in (\ref{local_operator}) by considering different values of $\delta$ defined in (\ref{delta}) and $p\equiv N_{sub}$. $\mathbf{u^{DA}}\in \mathbb{R}^{N_p\times N}$ is computed by running DD-4DVar algorithm for $N_{sub}=1$, while $\mathbf{\Tilde{u}^{DD-DA}}\in \mathbb{R}^{N_p\times N}$ is computed by gathering local solutions obtained by running DD-4DVAR algorithm for different values of $N_{sub}>1$ with $\delta\ge 0$, as shown in Figure (\ref{fig_DD_1D}). \\

\noindent In the following, we present experimental results of the consistency and stability of DD--4DVAR method considering the initial boundary problem of the SWEs 1D. The  discrete  model  is  obtained  using Lax--Wendroff scheme \cite{LeVeque} on $\Omega\times \Delta$ where the orders of convergence in space and time are $p=q=2$.
\\
\begin{itemize}
    \item Consistency of DD--4DVAR method.\\
    \\
\noindent 
\noindent From Table 1.1, we get\footnote{By fixing the values of $N_{sub}=N_t=4$ and $\delta=2$, we get $e^p(2,\Delta x,\Delta t, 4,4)={e}^{p}(\Delta x,\Delta t)$ i.e. $e^p$ defined in (\ref{error_exp}) depends only on $\Delta x$ and $\Delta t$.}
\begin{equation}\label{rel}
  e^{p}\left( \frac{\Delta x}{d}, \frac{\Delta t}{d}\right)\approx \frac{e^p(\Delta x,\Delta t)}{d^2}\quad d=1,2,4,6,8,10.
\end{equation}
As shown in Table 1.1 and Figure 1.6, \textit{the experimental order of  consistency corresponds to the theoretical one obtained in Theorem \ref{dd_truncation_error}.} \\
\item Stability of DD--4DVAR method.\\
\\
In Table 1.2 and Figure 1.7, we report values of $\bar{E}_k$ for different values of perturbation $\bar{e}_k$ on initial condition of $P_{i,k}^{M_{i,k}}$ defined in (\ref{co}). From the values in Table 1.2, we  experimentally estimate $C_k$ in (\ref{costante_C}), in particular 
\begin{equation*}
  C_k\approx 2.00\times 10^1 \quad \forall k=1,\ldots,N_t.
\end{equation*}
Consequently, \textit{DD--4DVAR with the initial boundary problem of SWEs 1D is well-conditioned.}

\end{itemize}

\begin{table}[ht!]
\centering
\begin{tabular}{ccccc}
\hline
$d$ & $\frac{\Delta x}{d}$ & $\frac{\Delta t}{d}$ & $e^p\left( \frac{\Delta x}{d}, \frac{\Delta t}{d}\right)$ & $\frac{e^p(\Delta x,\Delta t)}{d^2}$\\
 \hline
1&$7.87\times \times 10^{-3}$&$1.09\times 10^{-1}$ & $1.53\times 10^{-2}$ & $1.53\times 10^{-2}$\\
2&$3.92\times \times 10^{-3}$&$5.47\times 10^{-2}$ & $9.01\times 10^{-4}$ & $3.83\times 10^{-3}$\\
4&$1.96\times \times 10^{-3}$&$2.74\times 10^{-2}$ & $6.45\times 10^{-4}$ & $9.56\times 10^{-4}$\\
6&$1.30\times \times 10^{-3}$&$1.83\times 10^{-2}$ & $3.65\times 10^{-4}$ & $2.39\times 10^{-4}$\\
8&$9.78\times \times 10^{-4}$&$1.37\times 10^{-2}$ & $3.99\times 10^{-4}$ & $4.25\times 10^{-4}$\\
10&$7.81\times \times 10^{-4}$&$1.10\times 10^{-2}$ & $3.77\times 10^{-4}$ & $1.53\times 10^{-4}$\\
\hline
\end{tabular}
\caption{Fixed $N_p=640$ number of inner nodes in $\Omega$, $N=9$ number of instants of time in $\Delta$, $N_{sub}=4$ number of spatial subdomain and $N_t=4$ time intervals. We report the values of $e^p$ defined in (\ref{error_exp}) for different values of $\Delta x$ and $\Delta t$ spatial and temporal step sizes of $M_{i,k}$ defined in (\ref{local_matrix}).}
\label{table_value1}
\end{table}
\begin{table}[ht!]
\centering
\begin{tabular}{cc}
\hline
 $\bar{e}_k$&$\bar{E}_k$\\
 \hline
$3.03\times 10^{-6}$& $6.06\times 10^{-5}$\\
$3.03\times 10^{-5}$& $6.06\times 10^{-4}$\\
$3.03\times 10^{-4}$& $6.06\times 10^{-3}$\\
$3.03\times 10^{-3}$& $6.06\times 10^{-2}$\\
\hline
\end{tabular}
\caption{Fixed $N_p=640$ number of inner nodes in $\Omega$, $N=20$ number of instants of time in $\Delta$, $N_{sub}=4$ number of spatial subdomains and  $N_t=4$ time intervals. For $k=1,2,3,4$, we report the values of $\bar{E}_k$ defined in (\ref{def_prop_error}) for different values of perturbation $\bar{e}_k$ to initial condition of $P_{i,k}^{M_{i,k}}$ defined in (\ref{co}).}
\label{table_value7}
\end{table}

\begin{figure}
\centering
{\includegraphics[width=0.8\textwidth]{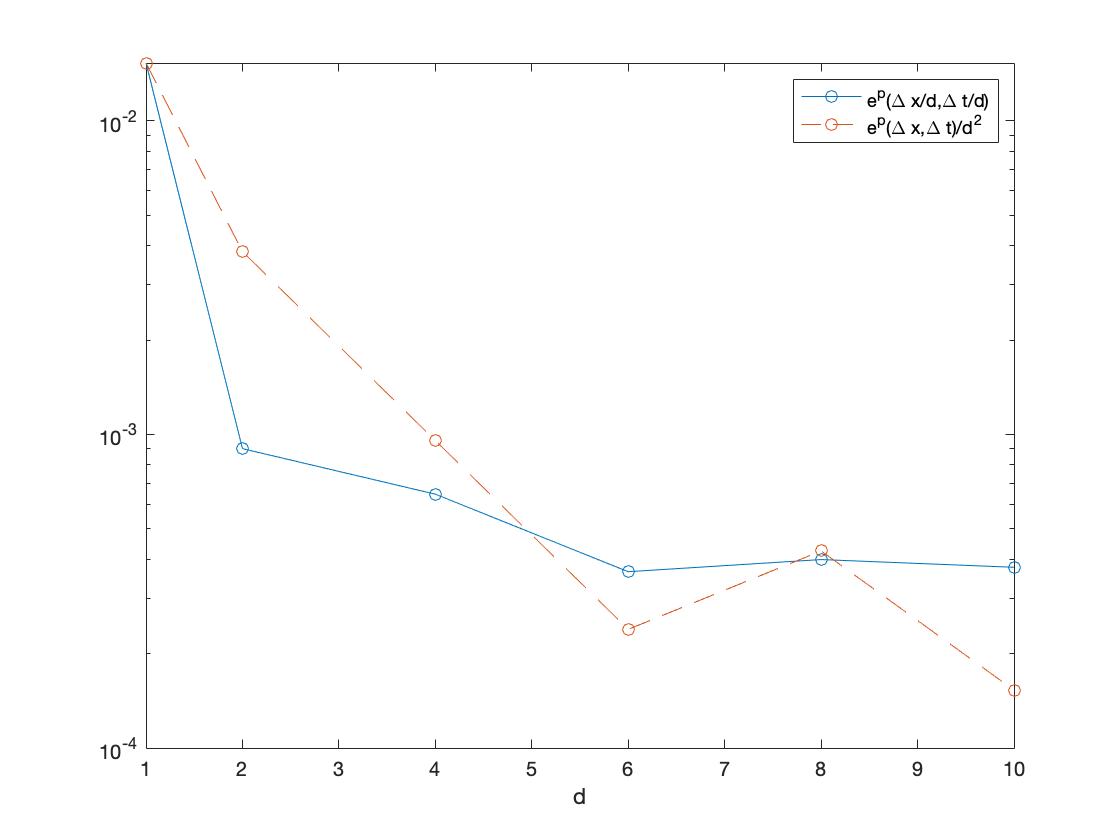}}\\
\caption{Plot of values of $e^p\left( \frac{\Delta x}{d}, \frac{\Delta t}{d}\right)$ (orange dashed line) and $\frac{e^p(\Delta x,\Delta t)}{d^2}$ (blue full line) for $d=1,2,4,6,8,10$ reported in Table 1.1.}
\label{plot_1}
\end{figure}

\begin{figure}
\centering
{\includegraphics[width=0.8\textwidth]{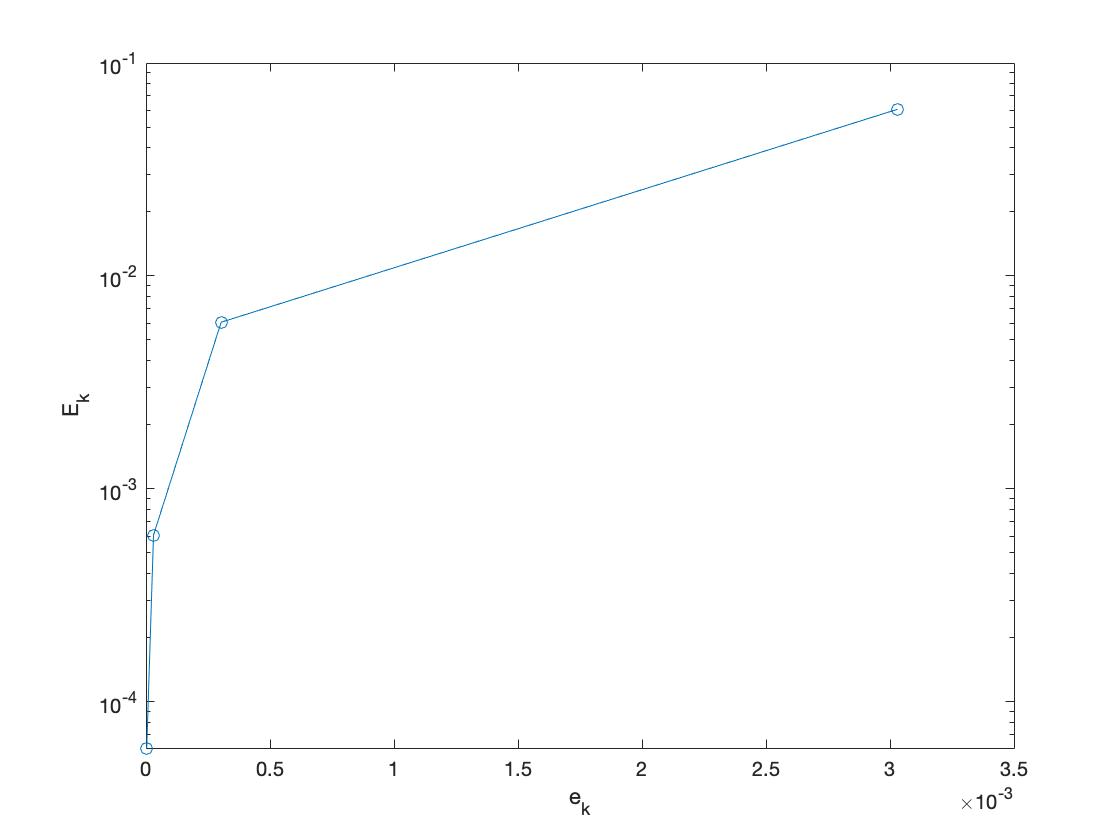}}\\
\caption{Plot of values $(\bar{e}_k,\bar{E}_k)$ reported in Table 1.2.}
\label{plot_2}
\end{figure}

\end{document}